\documentclass[12pt,a4paper]{article}

\usepackage[left=1cm,right=1cm,top=2cm,bottom=2cm]{geometry}
\usepackage[dvips]{graphics}
\usepackage[dvips]{epsfig}
\usepackage{color}
\usepackage{hyperref}
\usepackage{tikz,authblk}

\newtheorem{lemma}{Lemma}
\newtheorem{theorem}[lemma]{Theorem}
\newtheorem{corollary}[lemma]{Corollary}
\newtheorem{proposition}[lemma]{Proposition}
\newtheorem{definition}{Definition}
\newtheorem{remark}{Remark}
\newtheorem{conjecture}{Conjecture}
\newtheorem{example}{Example}

\newcommand{\dimo}{\noindent \emph{Proof. }}
\newcommand{\qed}{\\ \rightline{$\Box$ \ \ \ \ \ \ \ \ \ \ \ \ \ \ \ }\\}

\newcommand{\G}{\Gamma}
\newcommand{\g}{\gamma}

\usepackage{latexsym}
\usepackage{amsfonts}
\usepackage{amssymb}
\usepackage{amsmath}
\usepackage{graphicx}
\usepackage[english]{babel}

\begin{document}

\title{Crystallizations of compact 4-manifolds minimizing combinatorially defined PL-invariants}

 \renewcommand{\Authfont}{\scshape\small}
 \renewcommand{\Affilfont}{\itshape\small}
 \renewcommand{\Authand}{ and }

\author[1] {Maria Rita Casali}
\author[2] {Paola Cristofori}
\author[3] {Carlo Gagliardi}

\affil[1] {Department of Physics, Mathematics and Computer Science, University of Modena and Reggio Emilia, Via Campi 213 B, I-41125 Modena, Italy, casali@unimore.it}

\affil[2] {Department of Physics, Mathematics and Computer Science, University of Modena and Reggio Emilia, Via Campi 213 B, I-41125 Modena, Italy, paola.cristofori@unimore.it}

\affil[3] {Department of Physics, Mathematics and Computer Science, University of Modena and Reggio Emilia, Via Campi 213 B, I-41125 Modena, Italy, carlo.gagliardi@unimore.it}

\maketitle

\abstract {The present paper is devoted to present 
a unifying survey about some special classes of crystallizations 
of compact PL $4$-manifolds with empty or connected boundary,  called    \emph{ semi-simple} and \emph{ weak semi-simple crystallizations},  
with a particular attention to their properties of minimizing 
combinatorially defined PL-invariants,  such  
as the \emph{ regular genus}, the \emph{ Gurau degree}, the \emph{ gem-complexity} and the \emph{ (gem-induced) trisection genus}. \\
The main theorem, yielding a summarizing result on the topic, is an original contribution.  
Moreover, in the present paper the additivity of regular genus with respect to connected sum is proved to hold for all compact $4$-manifolds with empty or connected boundary 
which admit weak semi-simple crystallizations. }\endabstract

\bigskip
  \par \noindent
  {\small {\bf Keywords}: compact 4-manifold, colored triangulation, crystallization, regular genus, Gurau degree, gem-complexity, trisection}
  \smallskip
  \par \noindent
  {\small {\bf 2000 Mathematics Subject Classification}: 57Q15, 57N13, 57M15}

\section{Introduction} \label{intro}

It is well known that, thanks to a bright idea by Mario Pezzana (\cite{Pezzana_1}, \cite{Pezzana_2}),  every closed PL $n$-manifold $M$ can be triangulated by a pseudocomplex $K$, whose vertices are exactly $n+1$ (i.e. the minimum possible). 
If this is the case, $K$ and its dual edge-colored graph are called a {\it contracted triangulation} and a {\it crystallization} of $M$ respectively.

More recently,  the above result was extended to {\it singular $n$-manifolds}, 
i.e. triangulated polyhedra, whose vertices may have not only spheres, but also  closed connected  $(n-1)$-manifolds as links. In this context, some kind of ``minimality" with respect to the number of vertices of the obtained pseudocomplex can be considered, too. 
In particular, if $M$ has only one singular vertex, then Pezzana's theorem can be presented exactly in the same form. 
Hence, each such singular $n$-manifold may be combinatorially visualized and studied by means of regular graphs of degree $n+1$ (still called   {\it crystallizations}) whose edges are labelled by $n+1$ colors and such that the subgraph obtained by deleting all edges of any chosen color is connected  
(\cite{Cristofori-Mulazzani}, \cite{Casali-Cristofori-Grasselli}). 

Since singular $n$-manifolds with only one singular vertex are in bijection with manifolds with connected boundary, crystallizations can be thought of as a representation for manifolds with connected (non-spherical) boundary, too. 
Straightforward generalizations are known for singular $n$-manifolds with several singular vertices, i.e. for compact manifolds with several boundary components. 

\smallskip
The present paper is devoted to present a unifying survey about some special classes of crystallizations 
of compact PL $4$-manifolds with empty or connected boundary,   called    {\it semi-simple} and {\it weak semi-simple crystallizations}   (see Section 4 for details),   with a particular 
attention to their properties of minimizing interesting combinatorially defined PL-invariants,  
such  as the {\it regular genus}, the {\it Gurau degree} and the {\it gem-complexity}. 

\smallskip

The main achievement is the proof of the following summarizing result, which is an original contribution of the present paper.

\begin{theorem}[Main Theorem]\label{upper_bound_total}
Let $M^4$ be a compact $4$-manifold  with empty or connected boundary, 
and let $\widehat M^4$ be its associated singular manifold;    
let us assume 
\, $rk(\pi_1(M^4))= m \ge 0$ and $rk(\pi_1(\widehat M^4))= m^{\prime} \ge 0$. 
Then: 
\begin{itemize}
\item[(a)]
The regular genus $ \mathcal G (M^4)$ of $M^4$ satisfies     
$$ \mathcal G (M^4) \ge 2 \chi(\widehat M^4) + 5m - 2 (m- m^ {\prime}) -4.$$
Moreover, \  
equality holds if and only if $M^4$ admits a weak semi-simple crystallization. 
\item[(b)]
The Gurau-degree $ \mathcal D_G (M^4)$ of $M^4$ satisfies     
$$ \mathcal D_G (M^4) \ge 12 \Big[ 2 \chi(\widehat M^4) + 5m - 2 (m- m^ {\prime}) -4\Big]. $$
Moreover, \  
equality holds if and only if $M^4$ admits a semi-simple  crystallization.
\item[(c)]
The gem-complexity $k(M^4)$ of $M^4$ satisfies     
$$ k(M^4) \ge 3 \chi(\widehat M^4) + 10 m - 4 (m- m^ {\prime}) - 6.$$
Moreover, \  
equality holds if and only if $M^4$ admits a semi-simple crystallization.
\end{itemize}
\end{theorem}

\bigskip

A further original contribution of the paper is Proposition \ref{relation gem-complexity / regular genus}, yielding a characterization of compact $4$-manifolds admitting semi-simple crystallizations, via a relationship between gem-complexity and regular genus. 

\medskip 

In Section 5 the relevant\footnote{In the closed $4$-dimensional case, the problem is strictly related to the Smooth Poincar\'e Conjecture: see \cite{Ferri-Gagliardi} or Section 5.} 
problem of the additivity of regular genus with respect to connected sum is studied, and it is proved that the additivity holds for all compact $4$-manifolds with empty or connected boundary which admit weak semi-simple crystallizations: see Proposition \ref{prop:connected_sum}.  

\medskip 

Further, Section 6 recalls the notion of {\it gem-induced trisection} (due to \cite{Casali-Cristofori_trisections}), which extends the well-known notion of {\it trisection} (introduced in 2016 by Gay and Kirby: see \cite{Gay-Kirby}) to compact orientable $4$-manifolds with connected boundary, whose associated singular manifold is simply-connected. Also in this context, as a particular case of results proved in \cite{Casali-Cristofori_trisections}, weak semi-simple crystallizations turn out to have a ``minimality property", which enables to directly relate the so called {\it gem-induced trisection genus} with the  regular genus and/or the Betti numbers of the represented manifold:  see Propositions  \ref{trisection-weak} and \ref{minimality weak semi-simple}.

\bigskip

\section{Basic elements of crystallization theory} \label{preliminaries}

In the present section we will briefly review some basic notions of the so called {\it crystallization theory}, as a representation tool for 
piecewise linear (PL) compact manifolds; further details may be found in the quoted papers. 

From now on, unless otherwise stated, all spaces and maps will be considered in the PL category, and all manifolds will be assumed to be compact and connected.

\medskip

\begin{definition} \label{$n+1$-colored graph}
{\em An {\it $(n+1)$-colored graph}  ($n \ge 2$) is a pair $(\G,\g)$, where $\G=(V(\G), E(\G))$ is a multigraph (i.e. multiple edges are allowed, while loops are forbidden) 
which is regular of degree  $n+1$, and $\g$ is an {\it edge-coloration}, that is a map  $\g: E(\G) \rightarrow \Delta_n=\{0,\ldots, n\}$ which is injective on 
adjacent edges.}
\end{definition}

For sake of concision, when the coloration is clearly understood, colored graphs are often denoted simply by $\G$. 

\smallskip

For every  $\{c_1, \dots, c_h\} \subseteq\Delta_n$ let $\G_{\{c_1, \dots, c_h\}}$  be the subgraph obtained from $(\G, \g)$ by deleting all the edges that are not colored by the elements of $\{c_1, \dots, c_h\}$. 
In this setting, the complementary set of $\{c\}$ (resp. $\{c_1,\dots,c_h\}$)  in $\Delta_n$ will be denoted by $\hat c$ (resp. $\hat c_1\cdots\hat c_h$). 
The connected components of $\G_{\{c_1, \dots, c_h\}}$ are called {\it $\{c_1, \dots, c_h\}$-residues} or {\it $h$-residues} of $\G$; their number is denoted by $g_{\{c_1, \dots, c_h\}}$ (or, for short, by $g_{c_1,c_2}$, $g_{c_1,c_2,c_3}$ and $g_{\hat c}$ if $h=2,$ $h=3$ and $h = n$ respectively). 

 \medskip 

\noindent Each $(n+1)$-colored graph $\G$ encodes an associated $n$-dimensional pseudocomplex $K(\G)$:  
\begin{itemize}
\item $K(\G)$ contains an $n$-simplex for each vertex of $\G$, and the vertices of any $n$-simplex are (injectively) labelled by the elements of $\Delta_n$;
\item if two vertices of $\G$ are $c$-adjacent ($c\in\Delta_n$), then the corresponding $n$-simplices  of $K(\G)$ are glued along their $(n-1)$-dimensional faces opposite to the $c$-labelled vertices, so that equally labelled vertices are identified.
\end{itemize}

\smallskip

 In general $|K(\G)|$ is an {\it $n$-pseudomanifold} and $\G$ is said to {\it represent} it. 

\medskip

\noindent 
Via the above construction, it is not difficult to prove that: 
 \begin{itemize}
 \item[-] {\it $|K(\G)|$ is a closed $n$-manifold iff, for each color $c\in\Delta_n$, all $\hat c$-residues of $\G$ represent the $(n-1)$-sphere; 
 \item[-] $|K(\G)|$ is a singular
 \footnote{A polyhedron $|K|$ ($K$ being a simplicial complex) is said to be a {\it singular $n$-manifold} if the links of  the vertices of $K$ are closed connected $(n-1)$-manifolds.   The notion extends also to polyhedra associated to colored graphs: $|K(\Gamma)|$  is said to be a singular $n$-manifold if  the links of vertices of $K(\Gamma)$ in its first barycentric subdivision are closed connected $(n-1)$-manifolds. In both cases, a vertex whose link is not a $(n-1)$-sphere is called a {\it singular vertex}.} \ $n$-manifold iff, for each color $c\in\Delta_n$, all $\hat c$-residues of $\G$ represent closed connected $(n-1)$-manifolds.}
 \end{itemize}
\medskip

\begin{remark} \label{correspondence-sing-boundary} 
{\em Note that a bijective correspondence exists between singular $n$-manifolds and compact $n$-manifolds with no spherical boundary components. In fact, if $N$ is a singular $n$-manifold, then a compact $n$-manifold $\check N$ is easily obtained by deleting small open neighbourhoods of its singular vertices: $\check N$ turns out to be either closed (in case $N$ itself is a closed manifold, and hence  $N=\check N$) or with non-empty boundary, without spherical components.
Conversely, given a compact $n$-manifold $M$ without spherical boundary components, a singular $n$-manifold $\widehat M$ can be constructed by capping off each component of $\partial M$ by a cone over it.

\smallskip
For this reason, throughout the present work, we will restrict our attention to compact manifolds without spherical boundary components, and 
an $(n+1)$-colored graph $\G$ will be said to {\it represent} a compact $n$-manifold $M$ of this class (or, equivalently, to be a {\it gem} of $M$, where gem means {\it Graph Encoding Manifold}
: see \cite{Lins})  
if and only if  it represents the associated singular manifold $\widehat M$.}
\end{remark}

A restricted class of graphs gives the name to the whole theory: 

\begin{definition} {\em An $(n+1)$-colored graph $\Gamma$ representing a compact $n$-manifold with empty or connected boundary is said to be a {\it crystallization} of $M$ if,  for each color $c\in\Delta_n$, $\G_{\hat c}$ is connected.} 
\end{definition}


\medskip
The following theorem extends to the boundary case a well-known result - originally due to Pezzana (\cite{Pezzana_1}, \cite{Pezzana_2}) - founding the combinatorial representation theory for closed manifolds of arbitrary dimension via colored graphs. 

\begin{theorem}{\em 
(\cite{Casali-Cristofori-Grasselli}, \cite{Casali-Cristofori_generalized})}  
\ \label{Theorem_gem}  
Any compact orientable (resp. non orientable) $n$-manifold with no spherical boundary components admits a bipartite (resp. non-bipartite) $(n+1)$-colored graph representing it.
\ In particular, any  compact $n$-manifold with empty or connected boundary admits a crystallization representing it.   
\end{theorem}

\medskip

The  existence of a particular type of embedding of colored graphs into surfaces, is the key result in order to define two of the PL-invariants considered in the present paper.

\begin{proposition}{\em (\cite{Gagliardi_1981})}\label{reg_emb}
Let $\Gamma$ be a bipartite (resp. non-bipartite) $(n+1)$-colored graph of order $2p$. Then for each cyclic permutation $\varepsilon = (\varepsilon_0,\ldots,\varepsilon_n)$ of $\Delta_n$, up to inverse, there exists a cellular embedding of $\Gamma$  into an orientable (resp. non-orientable) closed surface $F_{\varepsilon}(\Gamma)$
whose regions are bounded by the images of the $\{\varepsilon_j,\varepsilon_{j+1}\}$-colored cycles, for each $j \in \mathbb Z_{n+1}$.
Moreover, the genus (resp. half the genus)  $\rho_{\varepsilon} (\Gamma)$ of $F_{\varepsilon}(\Gamma)$ satisfies

\begin{equation*}
2 - 2\rho_\varepsilon(\Gamma)= \sum_{j\in \mathbb{Z}_{n+1}} g_{\varepsilon_j,\varepsilon_{j+1}} + (1-n)p.   
\end{equation*}
\end{proposition}

\begin{definition} \label{Gurau-degree}
{\em Let $\Gamma$ be an $(n+1)$-colored graph.
If $\{\varepsilon^{(1)}, \varepsilon^{(2)}, \dots , \varepsilon^{(\frac {n!} 2)}\}$ is the set of all cyclic permutations of $\Delta_n$ (up to inverse), $ \rho_{\varepsilon^{(i)}}(\Gamma)$ \ ($i=1, \dots , \frac {n!} 2$) is called the {\it regular genus of $\Gamma$ with respect to the permutation $\varepsilon^{(i)}$}. Then, the {\it Gurau degree} (or {\it G-degree} for short) of $\Gamma$, denoted by  $\omega_{G}(\Gamma)$, is defined as
\begin{equation*}
 \omega_{G}(\Gamma) \ = \ \sum_{i=1}^{\frac {n!} 2} \rho_{\varepsilon^{(i)}}(\Gamma)
\end{equation*}
and the {\it regular genus} of $\Gamma$, denoted by $\rho(\Gamma)$, is defined as
\begin{equation*}
 \rho(\Gamma) \ = \ \min\, \{\rho_{\varepsilon^{(i)}}(\Gamma)\ /\ i=1,\ldots,\frac {n!} 2\}.
\end{equation*}}
\end{definition}

As a consequence, focusing on the represented compact $n$-manifolds, the following combinatorially defined PL-invariants are introduced: 

\begin{definition}\label{def_gen_degree} {\em Let $M$ be a compact (PL) $n$-manifold ($n\geq 2$).
The {\it (generalized) regular genus} of $M$ is defined as
\begin{equation*}
\mathcal G(M)=\min \{\rho(\Gamma)\ | \ \Gamma\mbox{ represents} \ M\}.   
\end{equation*}
and the {\it Gurau degree} (or {\it G-degree}) of $M$ is defined as
\begin{equation*}
\mathcal D_G(M)=\min \{\omega_G(\Gamma)\ | \ \Gamma\mbox{ represents} \ M\}.   
\end{equation*}}
\end{definition}

\begin{remark} \label{rem_regular genus} 
{\em Note that the (generalized) regular genus is a PL-invariant 
extending to higher dimension the classical genus of a surface and the Heegaard genus of a $3$-manifold. It succeeds in 
characterizing spheres in arbitrary dimension 
(\cite{Ferri-Gagliardi}), 
and a lot of classifying results via regular genus have been obtained, especially in dimension $4$ and $5$ (see \cite{Casali-Gagliardi_ProcAMS},  \cite{Casali_Forum2003}, \cite{Casali-Cristofori_generalized}  and their references). On the other hand, Gurau degree originally arises, within theoretical physics, from the theory of random tensors as an approach to quantum gravity in dimension greater than two (\cite{Gurau-book}).
Also G-degree characterizes spheres in arbitrary dimension and some classifying results via this invariant have recently been obtained in dimension $3$ and $4$:  see \cite{Casali-Cristofori-Grasselli} for the compact $3$-dimensional case, \cite{Casali-Cristofori-Dartois-Grasselli} for the closed $4$-dimensional case, and  \cite{Casali-Cristofori_generalized} for the compact $4$-dimensional case.}
\end{remark}

\medskip

A further  PL-invariant has been - quite naturally - defined within crystallization theory\footnote{Note that a lot of significant classification results have been obtained within crystallization theory with respect to gem-complexity, too:  as regards the closed case, see, for example, \cite{Bandieri-Cristofori-Gagliardi} and \cite{Casali-Cristofori_JKTR2008} for the dimension $3$, \cite{Casali-Cristofori_CATALOGUE} and \cite{Casali-Cristofori-Gagliardi_MONTESINOS} for the dimension $4$;  in the compact case, see \cite{Cristofori-Fomynikh-Mulazzani-Tarkaev} for a classification according to gem-complexity for compact orientable $3$-manifolds with toric boundary.}:

\par \noindent
\begin{definition} \label{gem-complexity}
{\em For each compact $n$-manifold $M$, its {\it gem-complexity} is the non-negative integer $k(M)= p - 1$, 
where $2p$ is the minimum order of an $(n+1)$-colored graph representing $M$.}  
\end{definition}

We point out that, for each compact $n$-manifold with empty or connected boundary, both regular genus and G-degree and gem-complexity are actually realized by 
a crystallization. 


Moreover,     
if  $M$ is a compact $n$-manifold with empty or connected boundary,  it is always possible to assume - up to a 
permutation of the color set - that any gem (and, in particular, any crystallization) of $M$ has color $n$ as its (unique) possible {\it singular color}, i.e. that each $\hat c$-residue, with $c \ne n$, represents the $(n-1)$-sphere.  

\medskip

In Section \ref{sec.trisections} a fourth PL-invariant (called {\it G-trisection genus}) will be combinatorially defined via colored graphs, in the restricted setting of compact $4$-manifolds $M^4$ such that the associated singular manifold $\widehat{ M^4}$  is simply-connected. 

\bigskip

\section{Computing invariants from crystallizations of compact 4-manifolds}

In the present section, $M^4$ will be a compact $4$-manifold with empty or connected  boundary, such that $rk(\pi_1(M^4))= m \ge 0$ and $rk(\pi_1(\widehat M^4))= m^{\prime} \ge 0$ (with $m^{\prime} \le m$), and  
$\Gamma$ will be a  $5$-colored graph representing $M^4.$   As pointed out in Section \ref{preliminaries}, we may assume without loss of generality $\Gamma$  to be a crystallization (i.e. $\Gamma_{\hat c}$ is connected 
for any $c \in \Delta_4$)  
and color $4$ to be its (unique) possible singular color (i.e. $\Gamma_{\hat c}$ represents $\mathbb S^3$,   
for any $c \ne 4$).  
Furthermore, let us denote by $\mathcal P_4$ the set of all cyclic permutations $\varepsilon=(\varepsilon_0,\varepsilon_1,\varepsilon_2,\varepsilon_3,\varepsilon_4)$ of $\Delta_4$ such that $\varepsilon_4=4.$

\medskip 

With the notations settled in Section \ref{preliminaries} for the number of residues,  \cite{Basak-Casali} and  \cite{Casali-Cristofori_generalized} yield, $ \forall j,k,l \in \Delta_3$: 
$$g_{j,k,l}=  1 + m^{\prime} + t_{j,k,l},  \ \ \text{with} \ t_{j,k,l} \ge 0  \ \ \text{and} \ \ \{r,s\}= \Delta_4 - \{j,k,l \};$$ 
$$g_{j,k,4}=  1 + m + t_{j,k,4},  \ \ \text{with} \ t_{j,k,4} \ge 0 \ \ \ \text{and} \ \ \{r,s\}= \Delta_3 - \{j,k\}.$$ 

As a consequence: 
\begin{equation} \label{sum g_ijk} 
\sum_{i,j,k \in \Delta_4}  g_{i,j,k} 
= 10 + 10m - 4 (m - m^{\prime}) + \sum_{i,j,k\in \mathbb{Z}_{5}}  t_{i,j,k}     
\end{equation}
\medskip

On the other hand,  in \cite{Casali-Cristofori_generalized} 
the following relation is proved to hold for each $i\in\Delta_4$ and for each $\varepsilon \in \mathcal P_4$: 

\begin{equation} \label{numerospigoli(n=4)}
g_{\widehat{\varepsilon_{i-1}},\widehat{\varepsilon_{i+1}}} = g_{\varepsilon_i,\varepsilon_{i+2},\varepsilon_{i+3}} \ =  \   1 
+     \rho_\varepsilon - \rho_{\varepsilon_{\widehat{i-1}}} - \rho_{\varepsilon_{\widehat{i+1}}}  
 \end{equation} 
 
\noindent where $\varepsilon_{\hat i} = (\varepsilon_0,\ldots,\varepsilon_{i-1},\varepsilon_{i+1},\ldots,\varepsilon_4=4)$ and  $\rho_\varepsilon,$   $\rho_{\varepsilon_{\hat{i}}}$ respectively denote   $\rho_\varepsilon(\G),$ $\rho_{\varepsilon_{\hat{i}}}(\G_{\widehat{\varepsilon_i}}). $ 
 
Therefore: 
$$ g_{\varepsilon_{i-1},\varepsilon_{i+1},\varepsilon_{i+3}} = 1 + \rho_\varepsilon  - \rho_{\varepsilon_{\hat i}}   -  \rho_{\varepsilon_{\widehat{i+2}}}  =  1+m^{\prime} + t_{\varepsilon_{i-1},\varepsilon_{i+1},\varepsilon_{i+3}}  \ \  \ \   \forall i \in \{2,4\} \ \ \text{and} $$
$$ g_{\varepsilon_{i-1},\varepsilon_{i+1},\varepsilon_{i+3}} = 1 + \rho_\varepsilon  - \rho_{\varepsilon_{\hat i}}   -  \rho_{\varepsilon_{\widehat{i+2}}}  =  1+m + t_{\varepsilon_{i-1},\varepsilon_{i+1},\varepsilon_{i+3}} \ \  \ \   \forall i \in \{0,1,3\},$$
which trivially imply 
\begin{equation} \label{genus-subgenera(t)}
\begin{aligned} \rho_{\varepsilon}  - \rho_{\varepsilon_{\hat i}}  -  \rho_{\varepsilon_{\widehat {i+2}}} -  m^{\prime} = & \ t_{\varepsilon_{i-1},\varepsilon_{i+1},\varepsilon_{{i+3}}}  \ \  \ \   \forall i \in \{2,4\}  \ \ \text{and} \\
 \rho_{\varepsilon}  - \rho_{\varepsilon_{\hat i}}  -  \rho_{\varepsilon_{\widehat {i+2}}} -  m =  & \ t_{\varepsilon_{i-1},\varepsilon_{i+1},\varepsilon_{{i+3}}}   \ \  \ \   \forall i \in \{0,1,3\} \end{aligned} \end{equation} 
  
\noindent where all subscripts are taken in $\mathbb Z_5.$   
\bigskip

Computations regarding the regular genus, the G-degree and the order of $\Gamma$, performed in the quoted papers and in \cite{Casali-Grasselli2019}, allow to prove the following summarizing result, which is an original contribution of the present paper.

\begin{proposition}\label{invariants computation}
Let $\Gamma$ be an order $2p$ crystallization of a compact $4$-manifold $M^4$ with empty or connected boundary, with \, $rk(\pi_1(M^4))= m \ge 0$ and $rk(\pi_1(\widehat M^4))= m^{\prime} \ge 0$ ($m^{\prime} \le m$). 
Then: 
\begin{itemize}
\item[(a)]
$$ \rho_{\varepsilon}(\Gamma) \  = \ 2 \chi(\widehat M^4) + 5m - 2 (m- m^ {\prime}) -4 + \sum_{i\in \mathbb{Z}_{5}}  t_{\varepsilon_i,\varepsilon_{i+2},\varepsilon_{i+4}};$$
\item[(b)]
$$ \omega_G(\Gamma)  \ = \ 6 \Big[ 4 \chi(\widehat M^4) + 10m - 4 (m- m^ {\prime}) -8 + \sum_{i,j,k\in \mathbb{Z}_{5}}  t_{i,j,k}\Big];$$
\item[(c)] 
$$ p-1  \ = \ 3 \chi(\widehat M^4) + 10m - 4 (m- m^ {\prime}) - 6 + \sum_{i,j,k\in \mathbb{Z}_{5}}  t_{i,j,k}.$$ 
\end{itemize}
\end{proposition}

\dimo
In \cite{Casali-Grasselli2019}, for each cyclic permutation $\varepsilon= (\varepsilon_0,\varepsilon_1,\varepsilon_2,\varepsilon_3,\varepsilon_4)$ of $\Delta_4$, the {\it associated permutation} $\varepsilon^{\prime}$ has been defined as $\varepsilon^{\prime} = (\varepsilon_0,\varepsilon_2,\varepsilon_4,\varepsilon_1,\varepsilon_3)$.\footnote{Note that, if  $\varepsilon \in \mathcal P_4$ is assumed (i.e. $\varepsilon_4=4$), we can always consider $\varepsilon^{\prime}= (\varepsilon_1,\varepsilon_3,\varepsilon_0,\varepsilon_2,\varepsilon_4=4)$, i.e. $\varepsilon^{\prime} \in \mathcal P_4$, too.} Then, \cite[Proposition 7]{Casali-Grasselli2019}  yields: 
\begin{equation}
\label{chi with associated permutations}
\chi(N^4)  \ = \ \Big(\rho_{\varepsilon}(\Gamma) + \rho_{\varepsilon^{\prime}}(\Gamma)\Big) -p + 3 
\end{equation}
for any  order $2p$ crystallization of a singular 4-manifold $N^4$ with one singular vertex at most. 

Moreover, in virtue of \cite[Proposition 6(b)]{Casali-Grasselli2019},  
$$\rho_{\varepsilon^{\prime}}(\Gamma) - \rho_{\varepsilon}(\Gamma)  \ = \   \sum_{j\in \mathbb{Z}_{5}} g_{\varepsilon_j,\varepsilon_{j+1},\varepsilon_{j+2}} - \sum_{j\in \mathbb{Z}_{5}} g_{\varepsilon_j, \varepsilon_{j+2}, \varepsilon_{j+4}}$$
holds for any $5$-colored graph representing a singular $4$-manifold $N^4$; hence:

\begin{equation}
\label{eq_difference}
\rho_{\varepsilon^{\prime}}(\Gamma) - \rho_{\varepsilon}(\Gamma) = \sum_{i\in \mathbb{Z}_{5}}  t_{\varepsilon_i,\varepsilon_{i+1},\varepsilon_{i+2}} - \sum_{i\in \mathbb{Z}_{5}}  
t_{\varepsilon_i,\varepsilon_{i+2},\varepsilon_{i+4}}.  
\end{equation}

Then, by comparing relations \eqref{eq_difference} and  \eqref{chi with associated permutations}, the following formula follows:
\begin{equation} \label{rho&chi}
\chi(\widehat M^4)  \ = \ 2 \rho_{\varepsilon}(\Gamma) + 3 -p + \sum_{i\in \mathbb{Z}_{5}}  t_{\varepsilon_i,\varepsilon_{i+1},\varepsilon_{i+2}} - \sum_{i\in \mathbb{Z}_{5}}   t_{\varepsilon_i,\varepsilon_{i+2},\varepsilon_{i+4}}    
\end{equation}

On the other hand, an easy computation (making use of \cite[Lemma 21]{Casali-Cristofori-Dartois-Grasselli}) yields: 
\begin{equation}
\label{chi with p}
\chi(\widehat M^4)= 5 - \frac 1 3  \sum_{i,j,k \in \Delta_4}  g_{i,j,k} + \frac 1 3 p.       
\end{equation}
Hence, by comparison with \eqref{rho&chi} and by using \eqref{sum g_ijk}: 
$$  \chi(\widehat M^4) = 2 \rho_{\varepsilon}(\Gamma) + 3 -3 \chi(\widehat M^4) + 5 -10 m + 4 (m- m^ {\prime}) - 2 \sum_{i\in \mathbb{Z}_{5}}  t_{\varepsilon_i,\varepsilon_{i+2},\varepsilon_{i+4}},$$   
from which 
\begin{equation} \label{rho-corrected}
\rho_{\varepsilon}(\Gamma)  = 2 \chi(\widehat M^4) + 5m - 2 (m- m^ {\prime}) -4 + \sum_{i\in \mathbb{Z}_{5}}  t_{\varepsilon_i,\varepsilon_{i+2},\varepsilon_{i+4}}
\end{equation}
easily follows, as well as
\begin{equation} \label{rho-epsilon'} \rho_{\varepsilon^{\prime}}(\Gamma)  = 2 \chi(\widehat M^4) + 5m - 2 (m - m^ {\prime}) -4 + \sum_{i\in \mathbb{Z}_{5}}  t_{\varepsilon_i,\varepsilon_{i+1},\varepsilon_{i+2}}.
\end{equation}
This proves statement (a). 

\medskip

In virtue of \cite[Proposition 5]{Casali-Grasselli2019}, 
$$ \omega_{G}(\Gamma)  \ = \ 6 \Big(\rho_{\varepsilon}(\Gamma) + \rho_{\varepsilon^{\prime}}(\Gamma)\Big)$$
holds for each $5$-colored graph  $(\Gamma, \gamma),$  and for each pair ($\varepsilon, \varepsilon^{\prime})$ of associated cyclic permutations of $\Delta_4$. 
Hence, by summing relations \eqref{rho-corrected} and \eqref{rho-epsilon'},  statement (b) easily follows:   
$$ \begin{aligned} 
\omega_G(\Gamma) \ & = \ 6 \Big( 2 \rho_{\varepsilon}(\Gamma) + (\rho_{\varepsilon^{\prime}}(\Gamma) - \rho_{\varepsilon}(\Gamma)) \Big) = \\
 \ & = \ 6\Big( 4 \chi(\widehat M^4) + 10 m - 4 (m- m^ {\prime}) -8 + \sum_{j,k,l \in \Delta_4}  t_{j,k,l} \Big).     
 \end{aligned}    $$

Finally, in order to prove statement (c), it is sufficient to make use of relation \eqref{chi with p}, together with relation  \eqref{sum g_ijk}: 
$$ \begin{aligned} 
p-1 \ & = \  3 \chi(\widehat M^4) -16 +  \sum_{i,j,k \in \Delta_4}  g_{i,j,k} \ = \\.    
\ &  = \   3 \chi(\widehat M^4) -16 + 10 + 10m - 4 (m - m^{\prime}) + \sum_{j,k,l \in \Delta_4}  t_{j,k,l}.    
 \end{aligned}    $$
\ \qed

\noindent  The following statement, extending \cite[Corollary 24]{Casali-Cristofori-Dartois-Grasselli} to the connected boundary case, 
is a direct consequence of Proposition \ref{invariants computation} (b) and (c):  
 
\begin{corollary}\label{relation between invariants}
Let $M^4$ be a compact $4$-manifold $M^4$ with empty or connected boundary, Then: 
$$  \mathcal D_G (M^4)  = 6 \Big[ \chi(\widehat M^4) - 2 + k(M^4)\Big].$$
\end{corollary}
\ \qed

\bigskip

\section{Weak semi-simple crystallizations of compact 4-manifolds}

In \cite{Basak-Casali} and \cite{Basak} two particular types of crystallizations are introduced and studied, by generalizing the notion of {\it simple crystallizations} for closed simply-connected $4$-manifolds (see \cite{Basak-Spreer} 
and \cite{Casali-Cristofori-GagliardiJKTR2015}): they are proved to be ``minimal" with respect to regular genus, among all graphs representing the same closed $4$-manifold.

In  \cite{Casali-Cristofori_trisections} these definitions are extended to compact $4$-manifolds with empty or connected  boundary.  

\begin{definition} \label{def_weak semi-simple} {\em Let $M^4$ be a compact $4$-manifold, with empty or connected boundary.  A $5$-colored graph $\G$ representing $M^4$   
is called {\it semi-simple}
if \, $g_{j,k,l} = 1 + m^{\prime} \ \ \forall \ j,k,l \in \Delta_3$ \, and \, $g_{j,k,4} = 1 + m \ \ \forall \ j,k \in \Delta_3,$ \, where  \, $rk(\pi_1(M^4))= m \ge 0$ and $rk(\pi_1(\widehat M^4))= m^{\prime} \ge 0$ ($m^{\prime} \le m$).     

\noindent $\G$ is called {\em weak semi-simple} with respect to a permutation $\varepsilon \in \mathcal P_4$ if \, 
$g_{\varepsilon_i,\varepsilon_{i+2},\varepsilon_{i+4}}$
$= 1 + m \ \ \forall \ i \in \{0,2,4\}$ \, and \, $g_{\varepsilon_i,\varepsilon_{i+2},\varepsilon_{i+4}} = 1 + m^{\prime} \ \ \forall \ i \in \{1,3\}$ (where the additions in subscripts are intended in $\mathbb{Z}_{5}$).
}\end{definition}

\noindent We point out that, as a consequence of the above definition, if $\G$ is weak semi-simple, then \ $g_{\hat j}=1,\ \forall \ j\in\Delta_4$, i.e. $\G$ is a crystallization of $M^4.$ 

\noindent In case $m=0$ (and, hence, $m^\prime=0$, too),  semi-simple (resp. weak semi-simple) crystallizations are said to be {\em simple} (resp. {\em weak simple}).

\medskip

By making use of relations  \eqref{genus-subgenera(t)}, for all $i \in \Delta_4$, it is not difficult to prove the following characterization of weak semi-simple crystallizations: 

\begin{proposition} {\rm (\cite[Corollary 8]{Casali-Cristofori_trisections})}\label{weaksemisimple}
Let $\G$ be a crystallization of a compact $4$-manifold $M^4$ with empty or connected boundary, 
with \, $rk(\pi_1(M^4))= m \ge 0$, $rk(\pi_1(\widehat M^4))= m^{\prime} \ge 0$ ($m^{\prime} \le m$). 
Then $\G$ is weak semi-simple with respect to 
a cyclic permutation $\varepsilon \in \mathcal P_4$
if and only if \  
$$ \rho_{\varepsilon_{\hat i}} = \frac 1 2 ( \rho_{\varepsilon} - m) \ \ \forall i \in  \Delta_3 \ \ \ \ \text{and}  \ \ \ \ \rho_{\varepsilon_{\hat 4}} = \frac 1 2 ( \rho_{\varepsilon} - m)+  ( m- m^{\prime}).$$
\end{proposition}
\medskip

\begin{example} \label{examples}
{\em As concerns the closed case, $\mathbb S^4,\ \mathbb{CP}^{2},\ \mathbb{S}^{2} \times \mathbb{S}^{2}$ admit simple crystallizations, while $\mathbb S^1 \times \mathbb S^3,\ \mathbb S^1 \widetilde \times \mathbb S^3$ (the orientable and non-orientable sphere bundles over $\mathbb S^1$)
and $\mathbb{RP}^4$ admit semi-simple crystallizations.   
See Figures 1, 2, 3, 4, 5 respectively.
Moreover, in \cite{Basak-Spreer} a simple (order 134) crystallization of the $K3$-surface is produced.  

In the boundary case, examples of simple crystallizations of  $\mathbb S^2 \times \mathbb D^2$ and $\xi_2$ -  the $\mathbb D^2$-bundle over $\mathbb S^2$ with Euler number $2$ , whose  boundary is the lens space $L(2,1)$  - are constructed in \cite{Casali-Cristofori_generalized}: see Figures 6 and 7.  

In the same paper, semi-simple crystallizations of $\mathbb Y^4_h$ and  $\widetilde{ \mathbb Y}^4_h$, the genus $h$ orientable and non-orientable $4$-dimensional handlebodies, can be found (see Figures 8 and 9, where the orientable cases $h=1$ and $h=2$ are depicted), as well as
a weak simple (but not simple!) crystallization of $\xi_c$ ($c \in \mathbb Z^+ - \{1,2\}$), the $\mathbb D^2$-bundle over $\mathbb S^2$ with Euler number $c$ whose boundary is the lens space $L(c,1)$: see Figure 10.

\smallskip

Other examples of weak simple crystallizations  may be found in the existing catalogue of rigid dipole-free bipartite crystallizations of closed orientable $4$-manifolds, up to 20 vertices (see \cite{Casali-Cristofori_CATALOGUE}): in particular, all elements with order 16 turn out to be weak simple crystallizations of simply-connected manifolds, whose simple crystallizations appear with less than 16 vertices.}

\begin{figure}[h]
\centerline{\scalebox{0.12}{\includegraphics{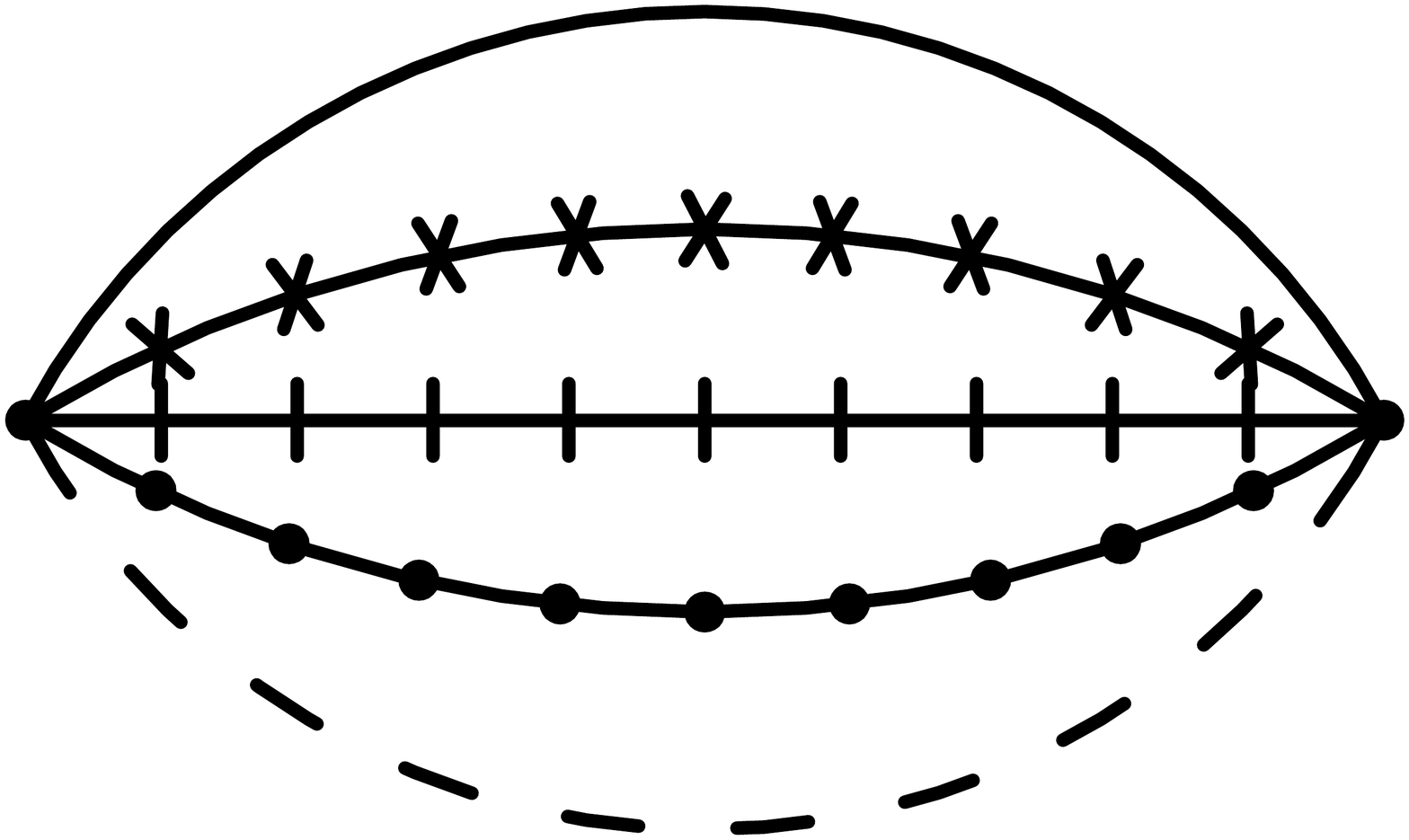}}}
\caption{\footnotesize The (unique) simple crystallization of  ${\mathbb S}^4$}
\label{fig.S4}
\end{figure}

\begin{figure}[h]
\centerline{\scalebox{0.6}{\includegraphics{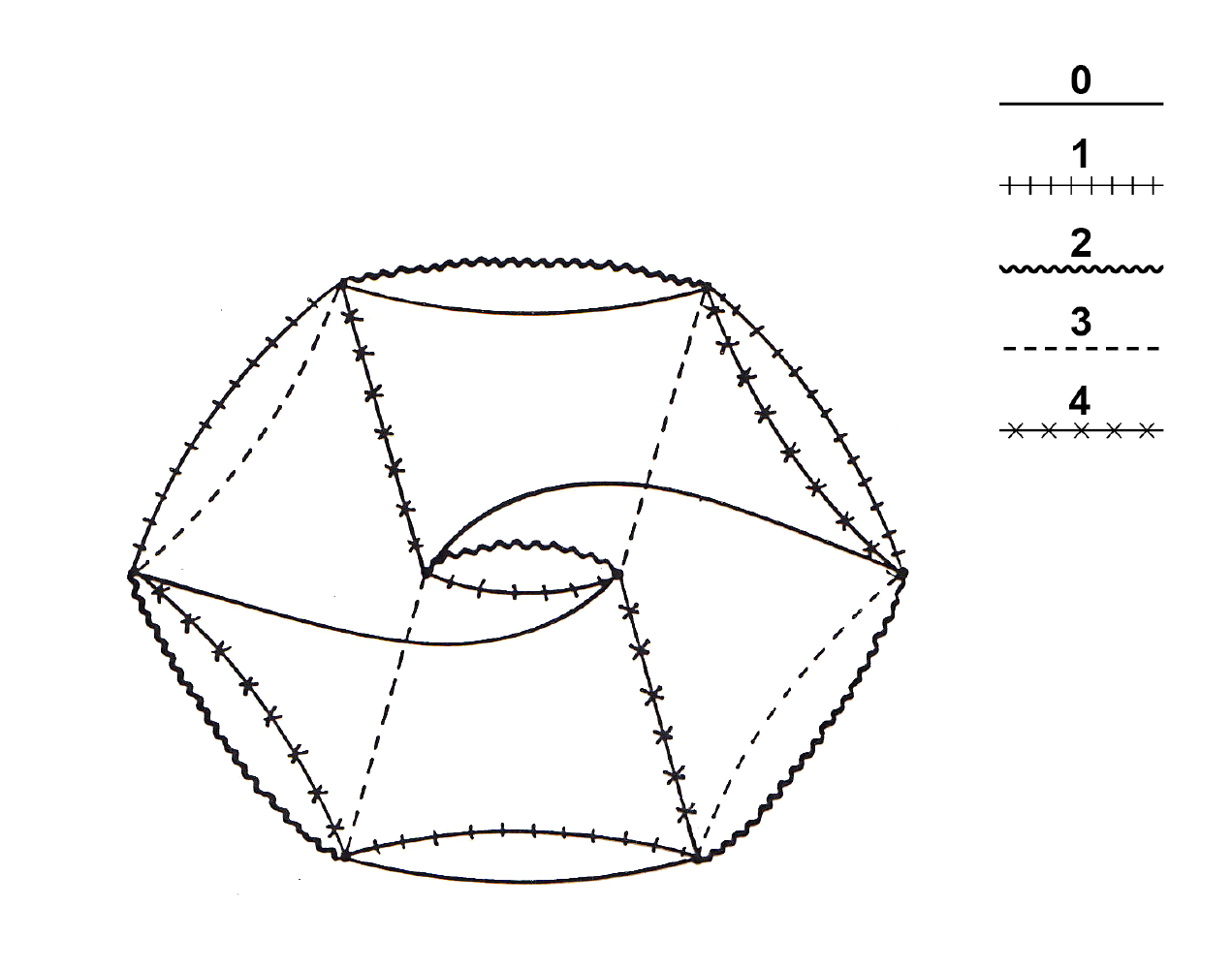}}}
\caption{\footnotesize The (unique) simple crystallization of  $\mathbb{ CP}^2$}
\label{fig.CP2}
\end{figure}

\begin{figure}[h]
\centerline{\scalebox{0.4}{\includegraphics{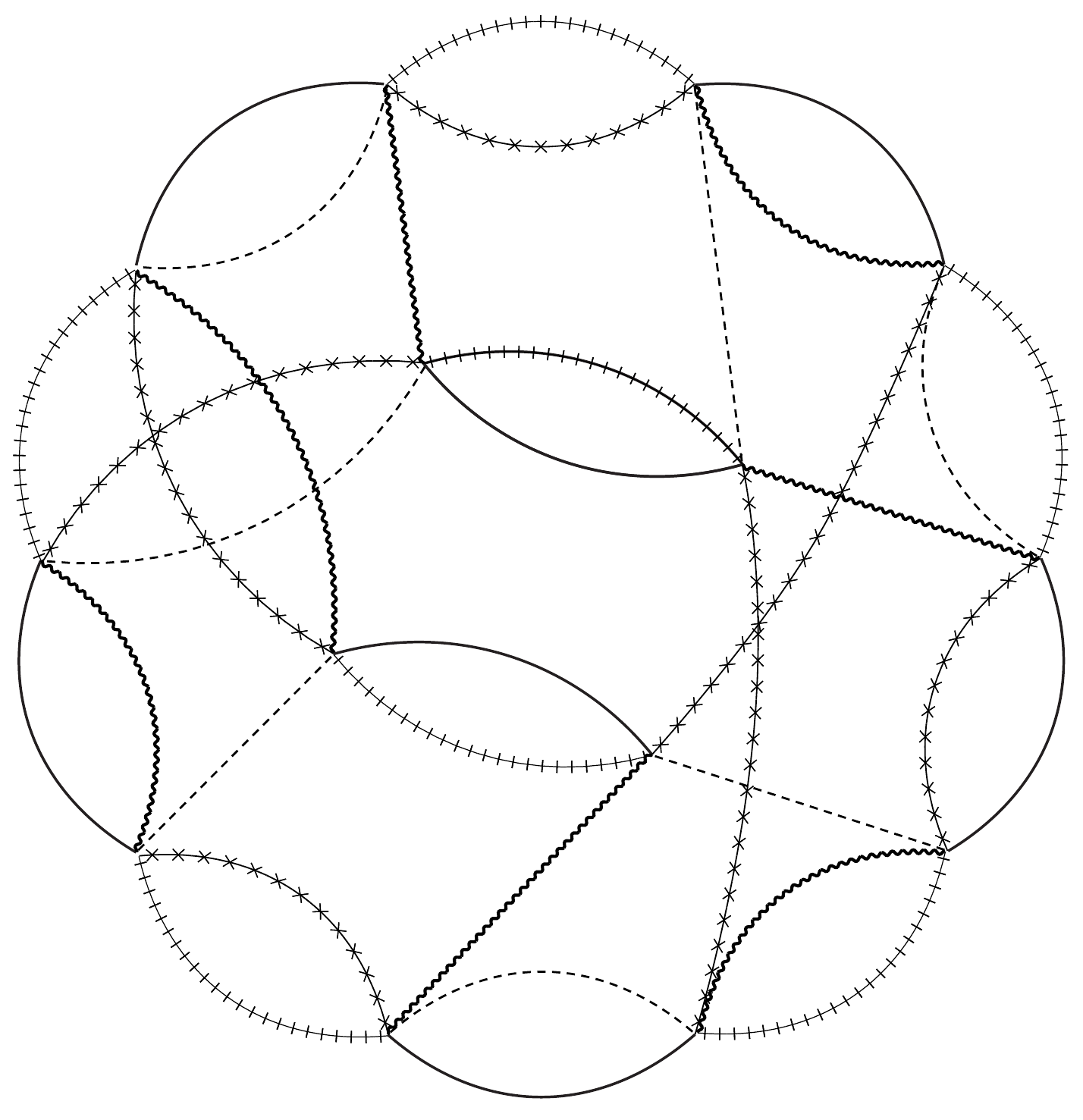}}}
\caption{\footnotesize A simple crystallization of  ${\mathbb S}^2 \times {\mathbb S}^2$}
\label{fig.S2xS2}
\end{figure}

\begin{figure}[h]
\centerline{\scalebox{1.1}{\includegraphics{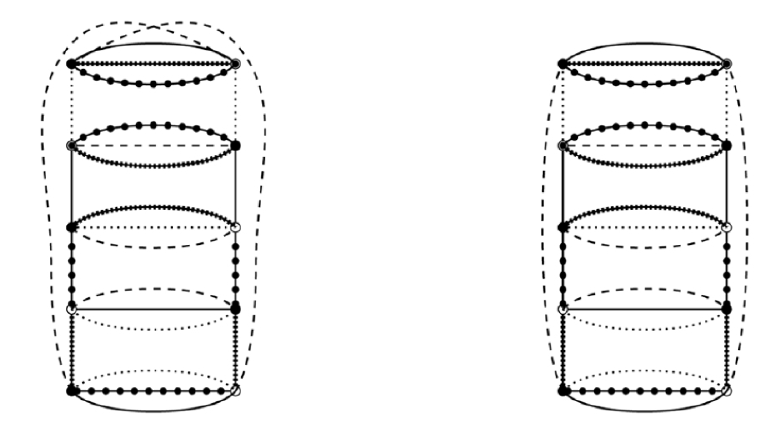}}}
\caption{\footnotesize Semi-simple crystallizations of  ${\mathbb S}^1 \times {\mathbb S}^3$ and  ${\mathbb S}^1 \ \widetilde \times \, {\mathbb S}^3$}
\label{fig.S3xS1_both}
\end{figure}

\begin{figure}[h]
\centerline{\scalebox{1.0}{\includegraphics{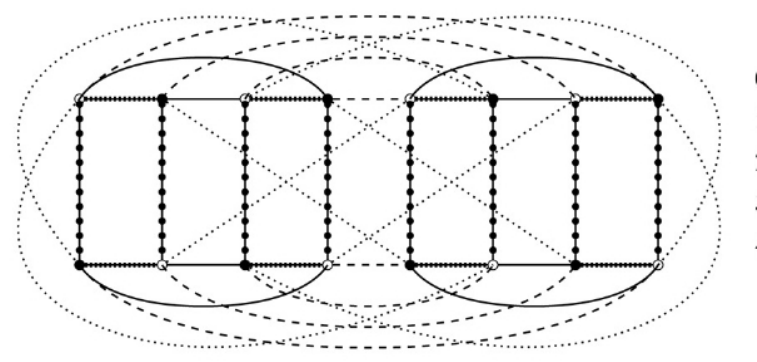}}}
\caption{\footnotesize The (unique) semi-simple crystallization of  $\mathbb {RP}^4$}
\label{fig.RP4}
\end{figure}

\begin{figure}[h]
\centerline{\scalebox{0.55}{\includegraphics{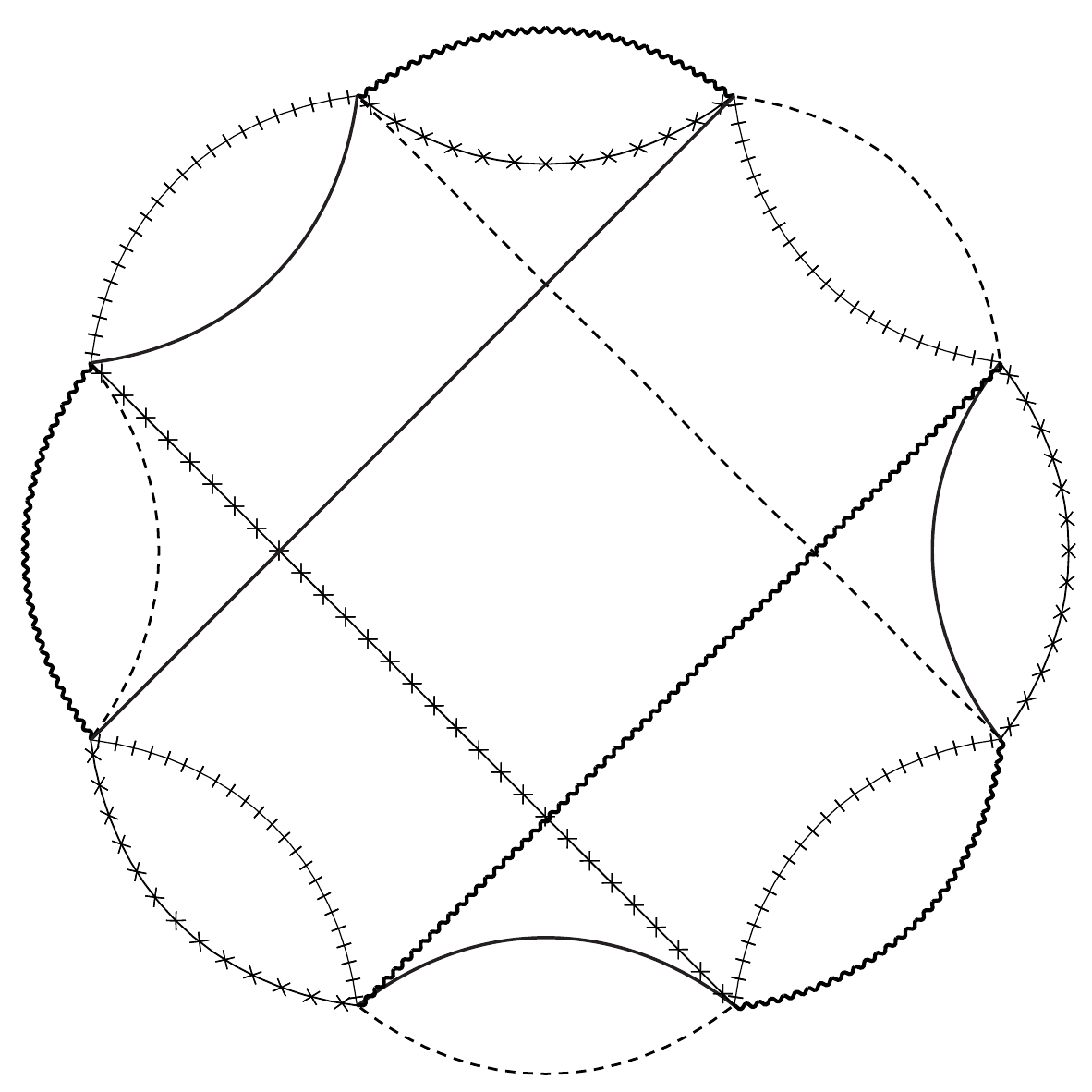}}}
\caption{\footnotesize A simple crystallization of  ${\mathbb S}^2 \times {\mathbb D}^2$ }  
\label{fig.S2xD2}
\end{figure}

\begin{figure}[h]
\centerline{\scalebox{0.55}{\includegraphics{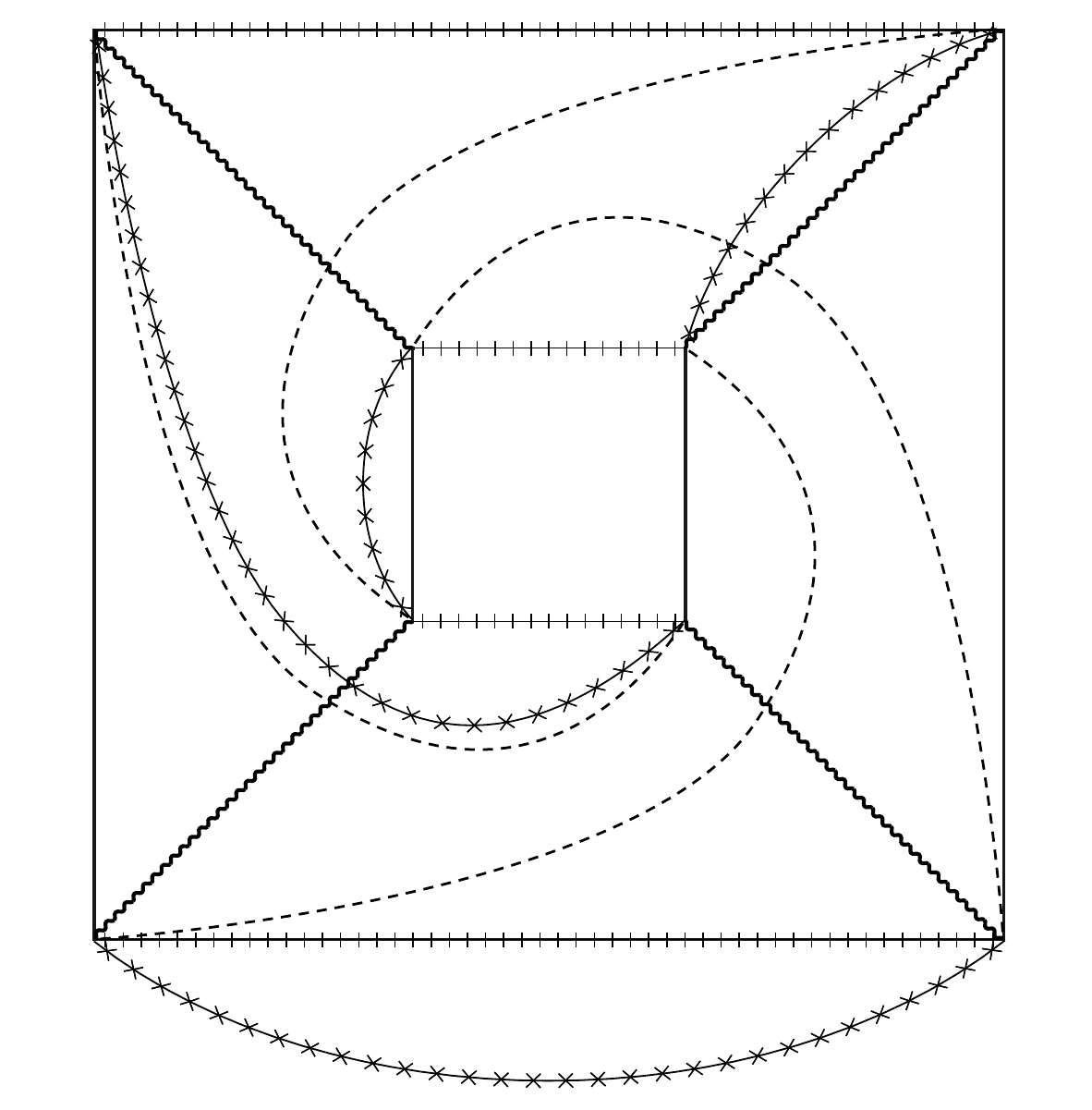}}}
\caption{\footnotesize A simple crystallization of  $ \xi_2$}
\label{fig.csi_2}
\end{figure}

\begin{figure}[h]
\centerline{\scalebox{0.55}{\includegraphics{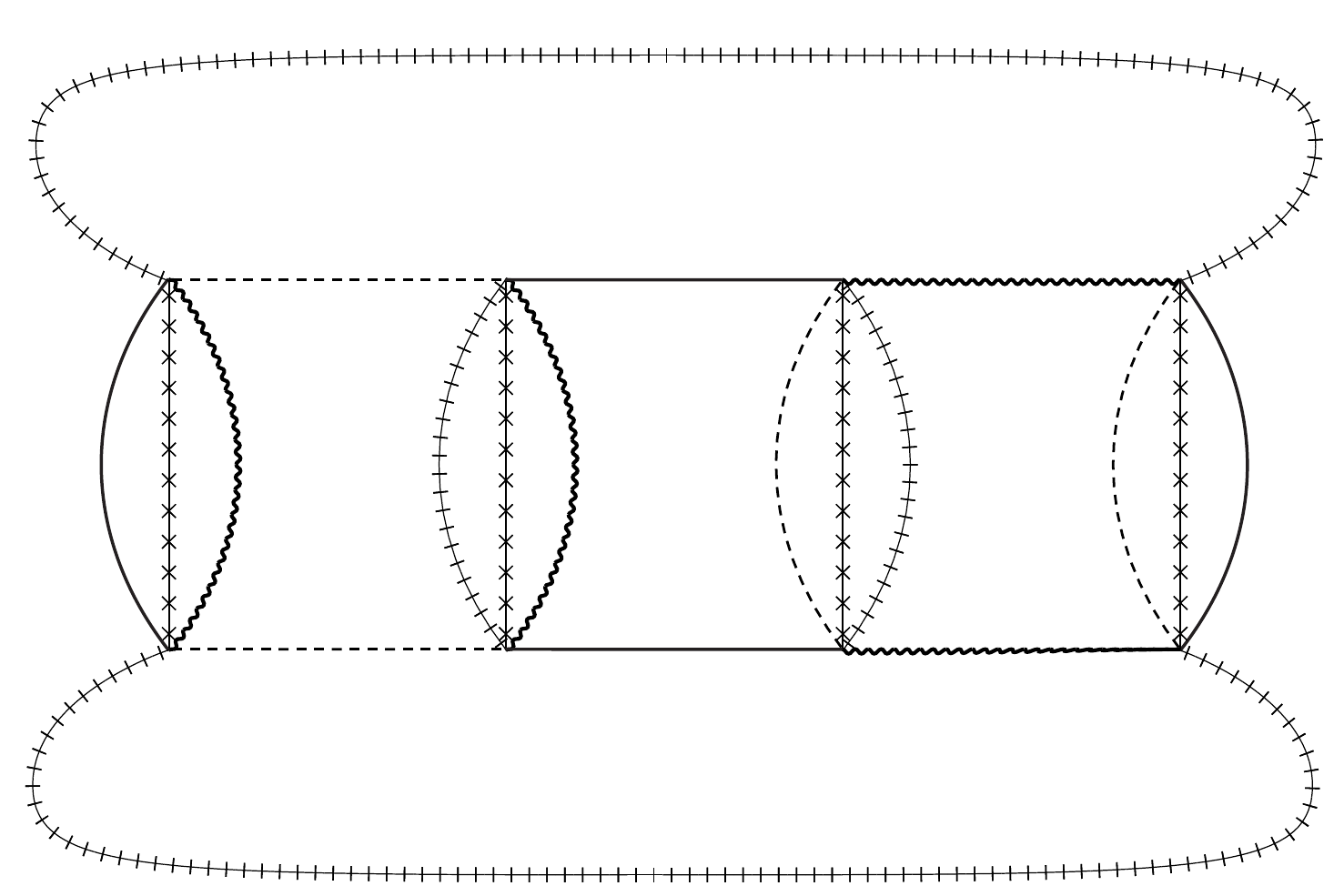}}}
\caption{\footnotesize A semi-simple crystallization of  $\mathbb Y^4_1$ }
\label{fig.Y^4_1}
\end{figure}

\begin{figure}[h]
\centerline{\scalebox{0.6}{\includegraphics{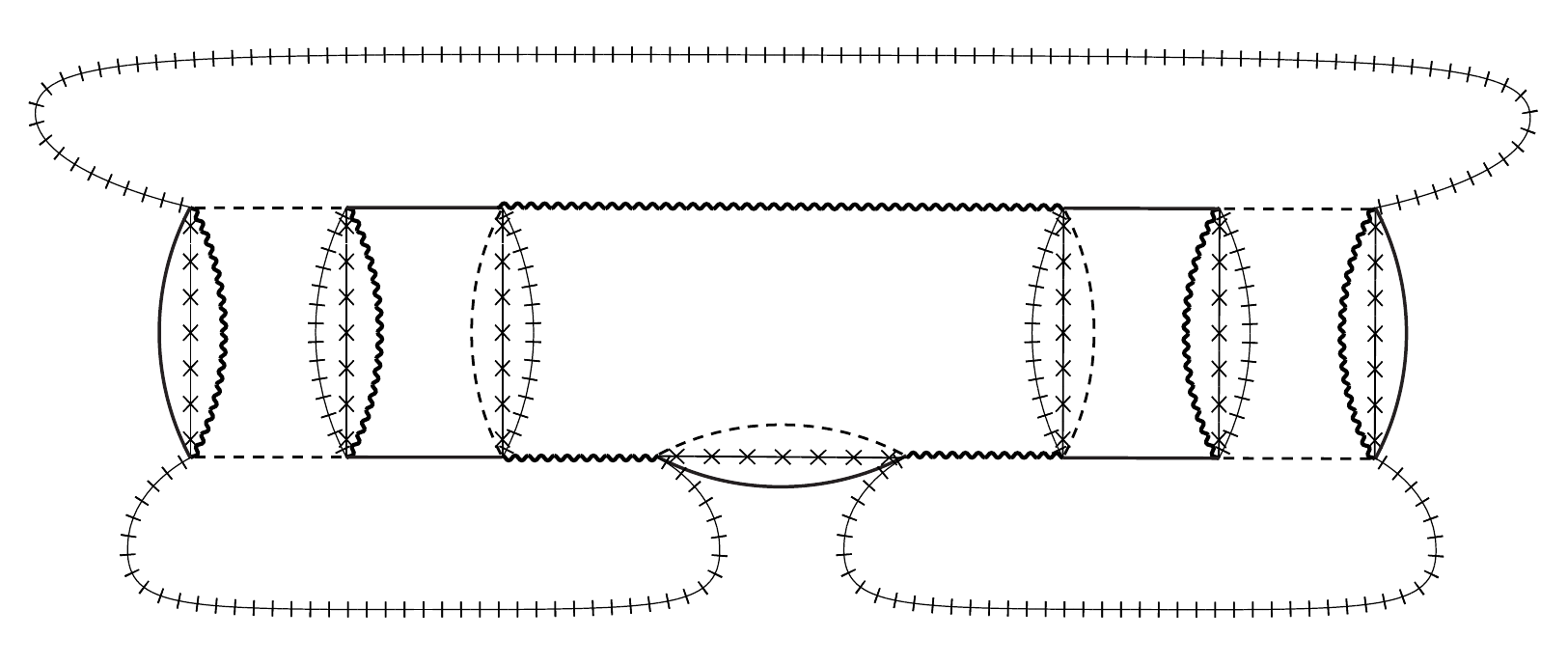}}}
\caption{\footnotesize A semi-simple crystallization of  $\mathbb Y^4_2$ }
\label{fig.Y^4_2}
\end{figure}

\begin{figure}[h]
\centerline{\scalebox{0.55}{\includegraphics{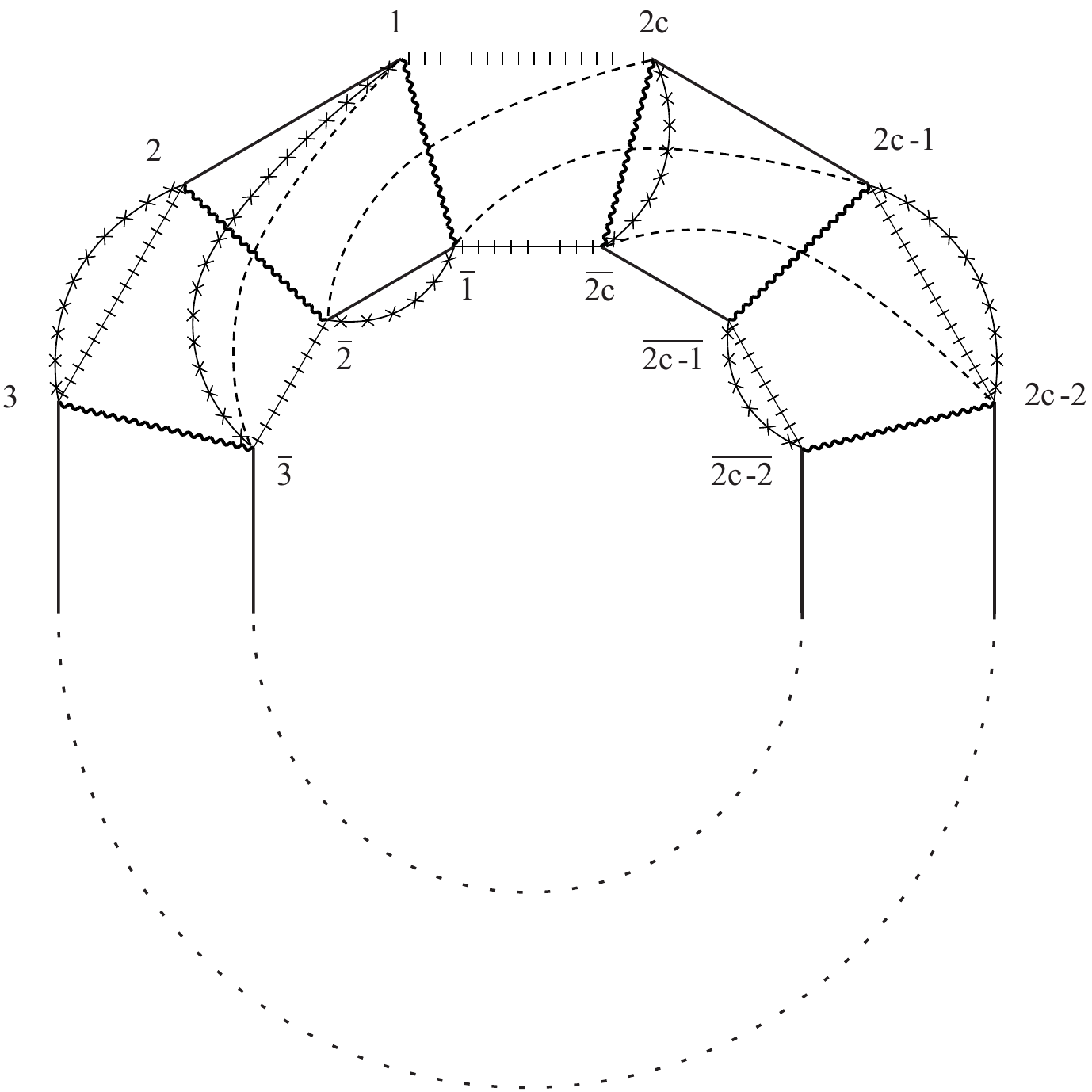}}}
\caption{\footnotesize A weak simple crystallization of  $ \xi_c$ ($c\geq 3$) }  
\label{fig.csi_c}
\end{figure}

\end{example}

We are now able to  prove the Main Theorem, stated in Section \ref{intro} . 
\medskip

\noindent {\it Proof of the Main Theorem.\ }
It is a direct consequence of  Proposition \ref{invariants computation}, together with the definitions themselves 
of semi-simple and weak semi-simple crystallization. 
In fact: 
$$\Gamma \ \text{weak semi-simple with respect to }  \  \varepsilon \in \mathcal P_4  \ \ \ \Longleftrightarrow  \ \ \ \ \sum_{i\in \mathbb{Z}_{5}}  t_{\varepsilon_i,\varepsilon_{i+2},\varepsilon_{i+4}}=0$$ 
$$\Gamma \ \text{semi-simple}  \  \ \ \ \Longleftrightarrow  \ \ \ \ \sum_{i,j,k \in \mathbb{Z}_{5}}  t_{i,j,k}=0$$ 
\vskip-0.4cm 
\ \qed

The following statement, characterizing manifolds which admit semi-simple crystallizations via a relationship between gem-complexity and regular genus, is an original contribution of the present paper. 

\begin{proposition}\label{relation gem-complexity / regular genus}
Let $M^4$ be a compact $4$-manifold with empty or connected boundary, with \, $rk(\pi_1(M^4))= m \ge 0$, $rk(\pi_1(\widehat M^4))= m^{\prime} \ge 0$ ($m^{\prime} \le m$).
Then: 
$$ k(M^4) = \frac{3 \mathcal G(M^4)+5m -2(m-m^{\prime})} 2 \ \    \Longleftrightarrow \ \  M^4  \ \text{admits a semi-simple 
crystallization}.$$
\end{proposition}

\dimo
Let $\Gamma$ and $\Gamma^\prime$ be two crystallizations of $M^4$ and $\varepsilon \in \mathcal P_4$ a permutation, such that  $\mathcal G (M^4) =  \rho (\Gamma) = \rho_\varepsilon (\Gamma)$ and  $k(M^4) = p^\prime -1$, $2p^\prime$ being the order of $\Gamma^\prime$. 
Statements (a) and (c) of Proposition \ref{invariants computation} yield:      
$$  \begin{aligned} \mathcal G (M^4) & = \rho_{\varepsilon}(\Gamma) \  = \ 2 \chi(\widehat M^4) + 5m - 2 (m- m^ {\prime}) -4 + \sum_{i\in \mathbb{Z}_{5}}  t_{\varepsilon_i,\varepsilon_{i+2},\varepsilon_{i+4}}, \\ 
k(M^4) & =  p^\prime -1  \ = \ 3 \chi(\widehat M^4) + 10m - 4 (m- m^ {\prime}) - 6 + \sum_{i,j,k\in \mathbb{Z}_{5}}  t^\prime_{i,j,k}, 
\end{aligned} $$
where $t_{i,j,k} \ge 0$ (resp. $t^\prime_{i,j,k} \ge 0$) is the difference between the number of $\{i,j,k\}$-residues in $\Gamma$ (resp. in $\Gamma^\prime$) and either $m+1$ (in case $4 \in  \{i,j,k\}$) or $m^\prime+1$ (in case $4 \notin  \{i,j,k\}$).    

Moreover, 
$ \sum_{i\in \mathbb{Z}_{5}}  t_{\varepsilon_i,\varepsilon_{i+2},\varepsilon_{i+4}} \le \sum_{i\in \mathbb{Z}_{5}}  t^\prime_{\bar\varepsilon_i,\bar\varepsilon_{i+2},\bar\varepsilon_{i+4}} \le  \sum_{i\in \mathbb{Z}_{5}}  t^\prime_{\bar\varepsilon_i,\bar\varepsilon_{i+1},\bar\varepsilon_{i+2}} $ 
holds, for any permutation $\bar\varepsilon\in \mathcal P_4$ such that $\rho_{\bar\varepsilon}(\Gamma^\prime)\leq\rho_{\bar\varepsilon^\prime}(\Gamma^\prime)$, $\bar\varepsilon^\prime$ denoting, as in the  of Proposition \ref{invariants computation}, the permutation of $\mathcal P_4$ which is associated to $\bar\varepsilon.$
\ Hence: 
$$ \begin{aligned} 2 k(M^4) & -  3  \mathcal G (M^4) - 5m + 2 (m- m^ {\prime})    \ =  \   2 \sum_{i,j,k\in \mathbb{Z}_{5}}  t^\prime_{i,j,k} - 3 \sum_{i\in \mathbb{Z}_{5}}  t_{\varepsilon_i,\varepsilon_{i+2},\varepsilon_{i+4}} \ge  \\
\ & \ge \  2 \sum_{i,j,k\in \mathbb{Z}_{5}}  t^\prime_{i,j,k}   - 3 \sum_{i\in \mathbb{Z}_{5}}  t^\prime_{\bar\varepsilon_i,\bar\varepsilon_{i+2},\bar\varepsilon_{i+4}} = \\
\ & = \    \sum_{i\in \mathbb{Z}_{5}}  t^\prime_{\bar\varepsilon_i,\bar\varepsilon_{i+1},\bar\varepsilon_{i+2}}  +  (\sum_{i\in \mathbb{Z}_{5}}  t^\prime_{\bar\varepsilon_i,\bar\varepsilon_{i+1},\bar\varepsilon_{i+2}}  - \sum_{i\in \mathbb{Z}_{5}}  t^\prime_{\bar\varepsilon_i,\bar\varepsilon_{i+2},\bar\varepsilon_{i+4}} ) \ \ge \ 0. 
\end{aligned}
$$
This proves that the equality $$2 k(M^4) =  3  \mathcal G (M^4) + 5m - 2 (m- m^ {\prime})$$ 
\noindent holds if and only if
$$\sum_{i\in \mathbb{Z}_{5}}  t^\prime_{\bar\varepsilon_i,\bar\varepsilon_{i+1},\bar\varepsilon_{i+2}} =0\quad\text{and}\quad\sum_{i\in \mathbb{Z}_{5}}  t^\prime_{\bar\varepsilon_i,\bar\varepsilon_{i+1},\bar\varepsilon_{i+2}}  - \sum_{i\in \mathbb{Z}_{5}}  t^\prime_{\bar\varepsilon_i,\bar\varepsilon_{i+2},\bar\varepsilon_{i+4}}=0$$ 
\noindent i.e. $\Gamma^\prime$ is a semi-simple crystallization of $M^4$.    
\ \qed

\medskip

As a consequence of  statement (a) of the Main Theorem, together with Proposition \ref{weaksemisimple}, we have: 

\begin{corollary}\label{weaksemisimple minimal genus}
Let $\G$ be a crystallization of a compact $4$-manifold $M^4$ with empty or connected boundary, with \, $rk(\pi_1(M^4))= m \ge 0$ and $rk(\pi_1(\widehat M^4))= m^{\prime} \ge 0$ ($m^{\prime} \le m$). 
Then, $\G$ is weak semi-simple with respect to the cyclic permutation $\varepsilon \in \mathcal P_4$ if and only if             
$$ \mathcal G (M^4) = \rho (\Gamma) = \rho_{\varepsilon}(\Gamma) = 2 \chi(\widehat M^4) + 5m - 2 (m- m^ {\prime}) -4$$
or, equivalently, if and only if 
$$ \rho_{\varepsilon_{\hat i}} = \chi(\widehat M^4) + m + m^ {\prime} -2 \ \ \forall i \in  \Delta_3 \ \ \ \ \text{and}  \ \ \ \ \rho_{\varepsilon_{\hat 4}} = \chi(\widehat M^4) + 2m -2.$$
\end{corollary}
\ \qed

\medskip

Table 1 summarizes the values of the invariants regular genus, G-degree and gem-complexity for all $4$-manifolds considered in Example \ref{examples}, computed via results of the present Section. 

\medskip

\begin{table}[!h]

\begin{center}
\begin{tabular}{|c||c|c|c|c|} 
\hline   \ & \ & \ & \ & \  \cr  
$M^4$ & $\mathcal G(M^4)$ & $\mathcal D_G(M^4)$ & $k(M^4)$  & notes \cr 
 \ & \ & \ & \ & \  \cr 
\hline \hline  
$\mathbb S^4$ & 0 & 0 &0 & admits simple crystallizations \cr  \hline
$\mathbb {CP}^2$ & 2 & 24 & 3 & admits simple crystallizations \cr  \hline
$\mathbb{S}^{2} \times \mathbb{S}^{2}$ & 4 & 48 &6 & admits simple crystallizations \cr  \hline
$\mathbb{S}^{1} \times \mathbb{S}^{3}$ & 1 & 12 &4 & admits semi-simple crystallizations \cr  \hline  
$\mathbb{S}^{1} \widetilde \times \mathbb{S}^{3}$ & 1 & 12 &4 & admits semi-simple crystallizations \cr  \hline  
$K3$ & 44 & 528 & 66 & admits simple crystallizations \cr  \hline  
$\mathbb{RP}^{4}$ & 3 & 36 & 7 & admits semi-simple crystallizations \cr  \hline
$\xi_2$ & 2 & 24 & 3 & admits simple crystallizations \cr  \hline
$\mathbb{S}^{2} \times \mathbb{D}^{2}$ & 2 & 24 & 3 & admits simple crystallizations \cr  \hline
$\mathbb Y^4_h$ & $h$ & 12$h$ & 3$h$ & admits semi-simple crystallizations \cr  \hline
$\widetilde{ \mathbb Y}^4_h$ & $h$ & 12$h$ & 3$h$ & admits semi-simple crystallizations \cr  \hline
$\xi_c$ ($c\in \mathbb Z^+-\{1,2\}$) & 2 & $\le 12 c$ & $\le 2c -1$ & admits weak simple crystallizations \cr  \hline
\end{tabular}
\caption{invariants for the considered compact $4$-manifolds}
\end{center}
\label{table_1}
\end{table}

\bigskip

We conclude the Section by pointing out what it is known about the topological structure of  simply-connected compact PL $4$-manifolds admitting simple or weak simple crystallizations. 
  
Recall that the notions of  simple and weak simple crystallizations arise from Definition \ref{def_weak semi-simple}, in case of $m=0$ (and, as a consequence, $m^\prime=0$, too).  

In particular, a simple crystallization (originally defined in \cite{Basak-Spreer} only in the closed case) is a $5$-colored graph representing a simply connected compact (PL) $4$-manifold with the property that the 1-skeleton of the associated triangulation equals the 1-skeleton of a $4$-simplex.

In \cite{Casali-Cristofori-GagliardiJKTR2015},  any (simply-connected) closed (PL) $4$-manifold $M$ admitting a simple crystallization is proved to admit a {\it special handlebody decomposition}, i.e. a handle decomposition lacking in 1-handles and 3-handles (see \cite[Section~3.3]{Mandelbaum}).  

\begin{theorem} {\rm \cite[Theorem 1.1]{Casali-Cristofori-GagliardiJKTR2015}}    \label{framed-link (closed)}
Let $M^4$ be a closed (PL) $4$-manifold. If $M^4$ admits a simple crystallization, then $M^4$ admits a handle decomposition lacking in 1-handles and 3-handles (or, equivalently, $M^4$ is represented by a (not dotted) framed link  with $\beta_2(M)$ components).
\end{theorem}

Recently,  the same property has been proved to hold also in the compact connected-boundary case, and also 
for a large class of $4$-manifolds, which comprehends those admitting  weak simple crystallizations.   

\begin{theorem} {\rm \cite{CC_handledec-weaksimple}}    \label{framed-link (weak)}
Let $M^4$ be a compact (PL) $4$-manifold, with empty or connected boundary. If $M^4$ admits a weak simple crystallization, then $M^4$ admits a handle decomposition lacking in 1-handles and 3-handles (or, equivalently, $M^4$ is represented by a (not dotted) framed link  with $\beta_2(M)$ components).
\end{theorem}

\begin{remark} {\em Note that - in the closed case - the existence of a special handlebody decomposition is related to Kirby problem n. 50: {\it ``Does every simply-connected closed $4$-manifold have a handlebody decomposition without 1-handles? Without 1- and 3-handles?"}. 
On the other hand, since simple crystallizations of TOP-homeomorphic PL manifolds must have the same order, the existence of infinitely many different PL-structures on the same TOP  $4$-manifold ensures that not all closed simply-connected PL $4$-manifolds admit simple crystallizations. Moreover, as a consequence of Theorem  \ref{framed-link (closed)}, it may be easily proved that, if an exotic PL-structure on $\mathbb S^4$ (resp. $\mathbb{CP}^{2}$) exists, then the corresponding PL $4$-manifold does not admit simple crystallizations (\cite[Corollary 3.3]{Casali-Cristofori-GagliardiJKTR2015}). 

Hence, Theorem \ref{framed-link (weak)} may be useful to investigate Kirby problem n. 50, via a class of 5-colored graphs
(including weak simple crystallizations) which possibly succeeds in representing all closed simply-connected $4$-manifolds. 
}
\end{remark}

\bigskip

\section{Invariants additivity  and related problems}
\label{sec.additivity}

It is well-known that the regular genus is subadditive with respect to connected sum of closed $n$-manifolds (\cite{FGG}).   
This can be checked directly via the so called {\it graph connected sum construction}, starting from a pair of graphs representing two given closed $n$-manifolds and realizing their 
regular genera. 
For each pair of $(n+1)$-colored graphs $\Gamma_1,$ $\Gamma_2$ and for each choice of vertices $v_1\in V(\Gamma_1)$, $v_2\in V(\Gamma_2)$, the {\it graph connected sum of $\Gamma_1$ and $\Gamma_2$ with respect to} 
$v_1$ and $v_2$ is the $(n+1)$-colored graph constructed by deleting $v_1$ and $v_2$ from $\Gamma_1$ and $\Gamma_2$ and welding the ``hanging'' edges of the same color.   
The obtained graph has regular genus equal to the sum of the regular genera of $\Gamma_1$ and $\Gamma_2$, and   - in case of $\Gamma_1,$ $\Gamma_2$ representing two closed $n$-manifolds $M_1$, $M_2$ - it is proved to represent a connected sum of $M_1$ and $M_2$.

\smallskip
On the other hand, the additivity of regular genus under connected sum has been conjectured\footnote{Obviously, in the 3-dimensional case, regular genus satisfies the additive property with respect to connected sum, via a 
classic result on Heegaard genus.}, and the associated (open) problem is significant, at least in the closed orientable case, and especially in dimension four.

\begin{conjecture} {\rm \cite{Ferri-Gagliardi}} \ \   
\label{I(n)}
Let $M^n_1, \ M^n_2$ be two closed (orientable) $n$-manifolds. \ \ Then, $$ \mathcal G(M^n_1  \# M^n_2) \, = \, \mathcal G(M^n_1) \, + \, \mathcal G(M^n_2).$$
\end{conjecture}

In fact, it is easy to prove that the 4-dimensional case of Conjecture \ref{I(n)} implies the 4-dimensional Smooth Poincar\'e Conjecture, via the well-known Wall Theorem on homotopic $4$-manifolds. 
   
\medskip 

However, it is obvious that the construction  of graph connected sum can be performed on any pair of $(n+1)$-colored graphs representing compact $n$-manifolds and, under suitable conditions, it yields an $(n+1)$-colored graph
representing either an internal or a boundary connected sum of the given manifolds (see \cite[Section 4]{Casali-Cristofori_generalized} for details).

Since the present paper focuses on compact manifolds with empty or connected boundary, in the following we will consider only the connected sum  constructions that are internal to this class of manifolds.

More precisely, given two compact $n$-manifolds $M^n_1,\ M^n_2$, with empty or connected boundary, $M^n_1\sharp M^n_2$ will denote an (internal) connected sum of $M^n_1$ and $M^n_2$ if and only if at least one of $M^n_1,\ M^n_2$ is closed; otherwise it will denote
a boundary connected sum of $M^n_1$ and $M^n_2$.  
Then, by exploiting the graph connected sum  construction, it is not difficult to check that 
$$\mathcal G(M^n_1 \sharp M^n_2) \le  \mathcal G(M^n_1) + \mathcal G(M^n_2).$$

\medskip

The definition itself of graph connected sum implies that the class of compact $4$-manifolds admitting weak simple/semi-simple crystallizations is closed under connected sum. This fact has interesting consequences regarding additivity properties of the PL-invariants regular genus, G-degree and gem-complexity. 

\begin{proposition} \ 
\label{prop:connected_sum}
Let $M^n_1$ and $M^n_2$ be two compact $4$-manifolds admitting weak semi-simple (resp. semi-simple) crystallizations.
Then,  $M^n_1 \sharp M^n_2$ admits weak semi-simple (resp. semi-simple) crystallizations, too.
\par \noindent
As a consequence, additivity of regular genus (resp. of G-degree and gem-complexity) holds within the class of compact $4$-manifolds admitting weak semi-simple (resp. semi-simple) crystallizations.

\ \qed \end{proposition}

Note that, in general,  the additivity of gem-complexity   (and of G-degree, too, via Corollary \ref{relation between invariants}) cannot hold because of the finiteness of the invariant: it is sufficient to make use of Wall Theorem, together with the existence of infinitely many PL-structures on the same TOP $4$-manifold.   

Notwithstanding this, Proposition \ref{prop:connected_sum} enables to compute both the regular genus and  the G-degree and the gem-complexity for a large class of compact $4$-manifolds, obtained by connected sums of the compact $4$-manifolds $\mathbb{CP}^{2}$, $\mathbb{S}^{2} \times \mathbb{S}^{2},$   $\mathbb{S}^{1} \times \mathbb{S}^{3}$, $\mathbb{S}^{1} \widetilde \times \mathbb{S}^{3},$  $K3,$ $\mathbb{RP}^{4}$, $\xi_2,$  $\mathbb{S}^{2} \times \mathbb{D}^{2},$ 
$\mathbb Y^4_h$, $ \widetilde{ \mathbb Y}^4_h$  (which admit semi-simple crystallizations, as shown in Example \ref{examples}).

\bigskip

The following statement extends to compact $4$-manifolds a double inequality concerning regular genus, obtained in \cite[Proposition 6.5]{Casali-Cristofori-Gagliardi_MONTESINOS}.

\begin{proposition} \label{extended double inequality}
For each compact $4$-manifold $M^4$ with empty or connected boundary, with $rk(\pi_1(M^4))=m \ge 0$ and  $rk(\pi_1(\widehat M^4))= m^{\prime} \ge 0$ ($m^{\prime} \le m$):
$$2 - 2 \mathcal G(M^4) \le \chi(\widehat M^4) \le 2 +  \frac  {\mathcal G(M^4)}2 - \frac {5m - 2(m - m ^{\prime})} 2.$$
\end{proposition}

\dimo
It is sufficient to make use of statement (a) of the Main Theorem, together with  the following formula, proved in  \cite[Proposition 13]{Casali-Cristofori_generalized}  for any 5-colored graph $\Gamma$ representing a  compact $4$-manifold $M^4$ and for each cyclic permutation $\varepsilon$ of $\Delta_4$: 
\begin{equation} \label{Euler-caracteristic n=4}
\chi (\widehat M^4)  = 2 - 2 \rho_{\varepsilon}(\Gamma) + \sum_{i \in \Delta_4} \rho_{\varepsilon_{\hat i}}(\Gamma_{\widehat{ \varepsilon_i}}).
\end{equation}
\vskip-0.3cm
\ \qed

\medskip

In \cite{Casali-Cristofori-Gagliardi_MONTESINOS}, by means of the double inequality improved by Proposition \ref{extended double inequality}, 
two classes of closed $4$-manifolds were detected, for which additivity of regular genus holds; 
while the complete identification of the manifolds belonging to one class was performed in the same paper, it was also pointed out that the problem of completely determining the other class remained open.
Now, we can trivially extend the analysis to the compact setting, obtaining the complete characterization of both (extended) classes. 

\begin{proposition} \label{additivity classes}
Let $M_1, M_2$ be two compact $4$-manifolds with empty or connected boundary, with $rk(\pi_1(M_i))=m_i\ \ge 0$ and $rk(\pi_1(\widehat M_i))= m_i^{\prime} \ge 0$ ($m_i^{\prime} \le m_i$)  for each $i \in \{1,2\}$.
\begin{itemize}
\item[(a)] If $\mathcal G(M_i) = 1 - \frac {\chi(\widehat M_i)} 2$ \ for each $i \in \{1,2\},$ \ then:
$$\mathcal G(M_1 \sharp M_2) =  \mathcal G(M_1) + \mathcal G(M_2) \ \ \   \text{and} \ \ \   \mathcal G(M_1\sharp M_2) = 1 - \frac {\chi(\widehat {M_1\sharp M_2})} 2.$$
\item[(b)] If $\mathcal G(M_i) = 2 \chi(M_i) + 5 m_i  -2(m_i - m_i^\prime) -4 $ \ for each $i \in \{1,2\},$ \ then:
$$\mathcal G(M_1 \sharp M_2) =  \mathcal G(M_1) + \mathcal G(M_2) \ \ \   \text{and}$$  
$$\mathcal G(M_1\sharp M_2) = 2 \chi(\widehat {M_1\sharp M_2}) + 5 (m_1 +m_2) -2(m_1 +m_2 - m_1^\prime -m_2^\prime)-4.$$   
\end{itemize}
Moreover, a compact $4$-manifold $M^4$ is involved in case (a) (resp. (b)) if and only if  it is a connected sum of sphere bundles over $\mathbb S^1$  
(resp. if and only if $M^4$ admits a weak semi-simple crystallization).
\end{proposition}

\dimo
Statements (a) and (b) are direct consequences of the double inequality  of Proposition \ref{extended double inequality} by further observing that $\chi(\widehat {M_1\sharp M_2}) = \chi(\widehat M_1) + \chi(\widehat M_2) -2.$ 
As regards the class of compact $4$-manifolds involved in statement (a), note that, in virtue of formula \eqref{Euler-caracteristic n=4}, they admit a 5-colored graph, realizing the regular genus (i.e. $\mathcal G(M^4)= \rho_{\varepsilon} (\Gamma)$), such that $\rho_{\varepsilon_{\hat i}} (\Gamma_{\widehat{ \varepsilon_i}})=0$ for each $i \in \Delta_4$. 
Now, \cite[Proposition 15]{Casali-Cristofori_generalized} ensures $M^4$ to be a  connected sum of (orientable or non-orientable) sphere bundles over $\mathbb S^1$.\footnote{Note that the condition $\rho_{\varepsilon_{\hat i}} (\Gamma_{\widehat{ \varepsilon_i}})=0$, $\forall i \in \Delta_4$ directly implies $M^4$ to be a closed $4$-manifold, 
since regular genus zero characterizes spheres in any dimension.} 

On the other hand, statement (a) of the Main Theorem easily proves that the class of compact $4$-manifolds involved in statement (b) exactly consists in compact $4$-manifolds admitting weak semi-simple crystallizations.
\ \qed

\smallskip
\begin{remark}
{\em We point out that the connected sums of (orientable or non-orientable) sphere bundles over $\mathbb S^1$ are the only compact $4$-manifolds belonging to both classes involved in Proposition \ref{additivity classes}: they are characterized by the equality between the regular genus and the rank of fundamental group  (as proved in \cite[Theorem 4]{Casali-Cristofori_generalized}), and hence  the double inequality of Proposition \ref{extended double inequality} actually becomes a double equality.} 
\end{remark}

\bigskip

\section{B-trisections induced by weak semi-simple crystallizations}
\label{sec.trisections}

Throughout this section all manifolds are supposed to be orientable.

The notion of {\it trisection} of a smooth, oriented closed $4$-manifold was introduced in 2016 by Gay and Kirby (\cite{Gay-Kirby}), by generalizing the classical 
idea of Heegaard splitting in dimension 3: the $4$-manifold is decomposed into three $4$-dimensional handlebodies, with disjoint interiors and mutually intersecting 
in $3$-dimensional handlebodies, so that the intersection of all three ``pieces" is a closed orientable surface. 
The minimum genus of the intersecting surface is called the {\it trisection genus} of the $4$-manifold. 

Hass, Bell, Rubinstein and Tillmann in \cite{Bell-et-al} performed an approach to the study of trisections via singular triangulations and their construction was applied by Spreer 
and Tillmann (\cite{Spreer-Tillmann}) to the case of triangulations induced by crystallizations of closed $4$-manifolds. In this setting Spreer and Tillmann succeeded into calculating the trisection genus
of all closed standard (PL) $4$-manifolds through their simple crystallizations.

The extension to the connected boundary case and to a wider class of edge-colored graphs is 
presented in \cite{Casali-Cristofori_trisections} following a suggestion in \cite{Rubinstein-Tillmann}; 
it relies on the notion of  B-trisection.

\begin{definition} \label{def_B-trisection} 
{\em Let $M^4$ be a compact $4$-manifold with empty (resp. connected)  boundary.  A  {\it B-trisection} of $M^4$  is a triple $\mathcal T =(H_{0},H_{1},H_{2})$ of $4$-dimensional submanifolds  of $M^4$, such that: 
\begin{itemize}
 \item [(i)]  $M^4 = H_{0}\cup H_{1}\cup H_{2}$ and $H_{0}, H_{1}, H_{2}$ have pairwise disjoint interiors;     
\item [(ii)]  $H_{1},H_{2}$ are $4$-dimensional handlebodies; $H_{0}$ is a $4$-disk \  (resp. is (PL) homeomorphic to $\partial M^4 \times [0,1]$);
 \item [(iii)] $H_{01}=H_{0}\cap H_{1},$  $H_{02}=H_{0}\cap H_{2}$   and $H_{12}=H_{1}\cap H_{2}$    are $3$-dimensional handlebodies;
\item [(iv)] 
$\Sigma (\mathcal T) = H_{0}\cap H_{1}\cap H_{2}$ is a closed connected surface (which is called {\it central surface}). 
\end{itemize}
}
\end{definition}
\smallskip

\begin{remark} \label{rem-simply-connected} 
{\em Note that the central surface of a B-trisection $\mathcal T = (H_{0},H_{1},H_{2})$ of $M^4$ is an Heegaard surface for the 3-manifold $\partial H_i=\#_{k_i}  (\mathbb S^1 \times \mathbb S^2)$ ($k_i \ge 0$), for each $i\in \{1,2\}$, splitting it into the 3-dimensional handlebodies $H_{ij}$ and $H_{ik}$, with $\{j,k\}= \{0,1,2\}-\{i\}$. 
Moreover, in the closed (resp. boundary) case, $(H_{01},H_{02},\Sigma (\mathcal T))$ is an Heegaard splitting of $\partial H_0 = \mathbb S^3$ (resp. of $\partial M^4$, and more precisely of the boundary component of $\partial H_0$ intersecting $H_1\cup H_2$).
\\\noindent
Hence, obviously, we have $ k_i \le genus(\Sigma (\mathcal T))$ for each $i\in \{1,2\}$, and, in the boundary case, the genus of $\Sigma (\mathcal T)$ is an upper bound for the Heegaard genus of $\partial M^4$. 
\\  \noindent
Moreover, via Seifert-Van Kampen's Theorem, it is not difficult to check that, the simply-connectedness of the singular manifold $\widehat{M^4}$ is a necessary condition for the 
existence of a B-trisection of $M^4$(see \cite{Casali-Cristofori_trisections}).}\end{remark} 

\medskip 

In order to construct B-trisections for  $4$-manifolds with empty or connected boundary, we consider the set, denoted by $G_s^{(4)}$, of all $5$-colored graphs having only one $\hat 4$-residue and such that all $\hat i$-residues, with $i\in\Delta_3$, represent the 3-sphere. 
Note that any compact $4$-manifold with empty or connected boundary  can be represented by an element of this set. 
Moreover,  $G_s^{(4)}$ properly contains (up to permutation of the color set)  all weak semi-simple crystallizations of compact $4$-manifolds with empty or connected boundary.  

\smallskip
The following theorem 
ensures the existence, for the whole class of compact $4$-manifolds with empty or connected boundary, 
of a triple of submanifolds satisfying ``almost all" conditions required by a B-trisection.

\begin{theorem} {\rm \cite{Casali-Cristofori_trisections}} \label{theorem_B-trisection}  
Let $M^4$ be a compact $4$-manifold with empty or connected boundary. 
For each 5-colored graph $(\Gamma,\gamma)\in G_s^{(4)}$ representing $M^4$ and for 
each $\varepsilon\in\mathcal P_4$, a  triple $\mathcal  T(\Gamma, \varepsilon) =(H_{0},H_{1},H_{2})$ of submanifolds of $M^4$ is obtained, satisfying properties (i), (ii) and (iv) of Definition \ref{def_B-trisection}, and such that  $H_{01}=H_{0}\cap H_{1}$ and 
 $H_{02}=H_{0}\cap H_{2}$  are $3$-dimensional handlebodies. 
\\ Moreover, the central surface $\Sigma (\mathcal T (\Gamma, \varepsilon)) = H_{0}\cap H_{1}\cap H_{2}$  is a closed connected surface of genus $\rho_{\varepsilon_{\hat 4}}(\Gamma_{\hat{\varepsilon_4}})$.    
\end{theorem}

The notions of gem-induced trisection and G-trisection genus arise quite naturally from the above result.    

\begin{definition} \label{def_gem-induced trisection} 
{\em Let $M^4$ be a compact $4$-manifold $M^4$ with empty or connected  boundary. If the triple $\mathcal T (\Gamma, \varepsilon)=(H_{0},H_{1},H_{2})$ of $M^4$, associated to a 5-colored graph $\G\in G_s^{(4)}$  and a permutation $\varepsilon \in \mathcal P_4,$ 
is a B-trisection (i.e. if  $H_{12}=H_{1}\cap H_{2}$ is a $3$-dimensional handlebody, too), then it is called  a {\it gem-induced trisection} of $M^4.$}

\smallskip
  
\noindent {\em The  {\it gem-induced trisection genus} -  or  {\it G-trisection genus} for short -  $g_{GT}(M^4)$ of $M^4$ is the minimum genus of the central surface of any gem-induced trisection $\mathcal T (\Gamma, \varepsilon)$ of $M^4$: 
 $$   g_{GT}(M^4) \ = \  \min \{genus(\Sigma (\mathcal T (\Gamma, \varepsilon)))  \ /  \ \mathcal T (\Gamma, \varepsilon)  \ \text{B-trisection of} \ M^4\}.$$ 
}\end{definition}
  
\smallskip

As a direct consequence of Theorem \ref{theorem_B-trisection}, if $\G$ is a crystallization of a closed simply-connected $4$-manifold $M^4$  
admitting a gem-induced trisection - actually a trisection -, then the trisection genus of $M^4$ is less or equal to $ \rho_{\varepsilon_{\hat 4}}(\Gamma_{\hat{\varepsilon_4}}).$

Detecting classes of $5$-colored graphs inducing triples $\mathcal T (\Gamma, \varepsilon)$, and possibly  B-trisections, with minimal genus of their central surface is thus a relevant problem.
The problem is faced in \cite{Casali-Cristofori_trisections}, proving in particular that weak semi-simple crystallizations guarantee the minimality of the genus of the central surface of a gem-induced trisection, if any.   
For this reason, we will now briefly sketch the construction of the triple $\mathcal  T(\Gamma, \varepsilon)$ of Theorem \ref{theorem_B-trisection} for the particular case of a crystallization $\Gamma$ that is weak semi-simple with respect to the permutation $\varepsilon\in\mathcal P_4.$

\medskip

Let us denote by $\sigma$ the standard $2$-simplex, by $v_b,v_r,v_b$ its vertices and by $\sigma^\prime$ its first barycentric subdivision.
Following \cite{Bell-et-al} and \cite{Spreer-Tillmann}, let us consider the following partition of $\Delta_4:$ $b=\{4\},\ g=\{\varepsilon_1,\varepsilon_3\},\ r=\{\varepsilon_0,\varepsilon_2\}.$ 
In the following, for sake of simplicity, we suppose $\varepsilon = (0,1,2,3,4)$.

\begin{figure}[t]
\centering
\scalebox{0.5}{\includegraphics{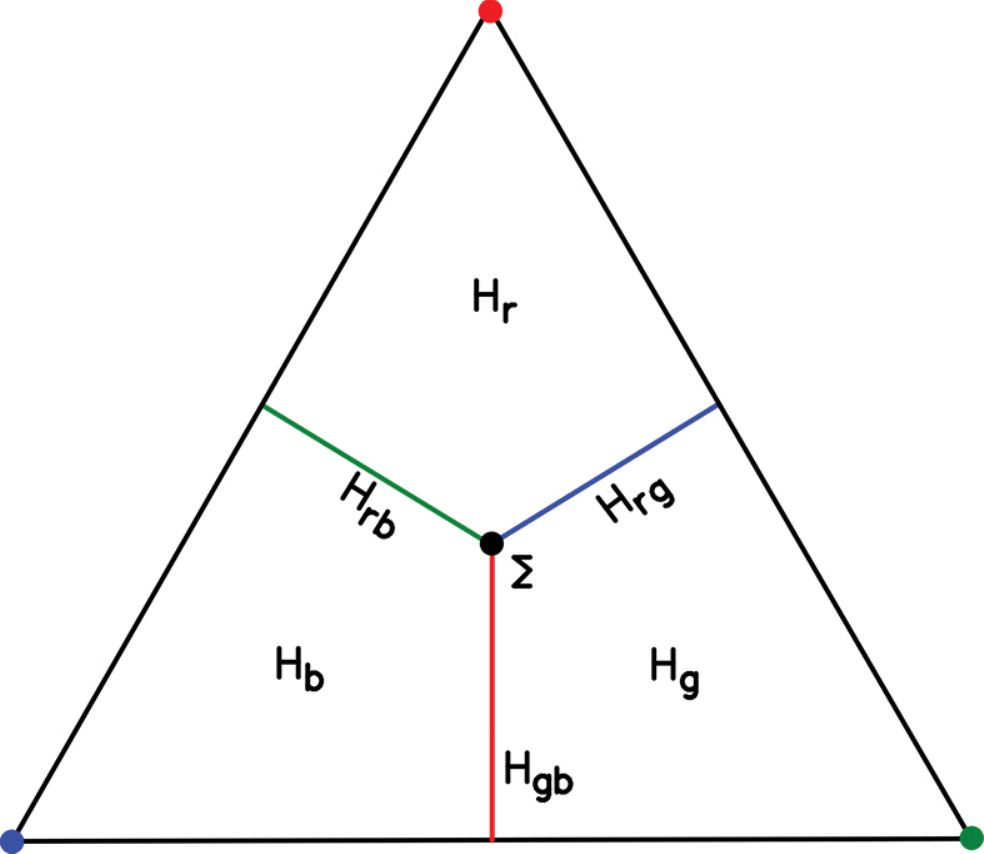}}
\caption{The 2-simplex $\sigma$}
\label{figsigma}
\end{figure}

Let $\mu\ :\ K(\G)\to\sigma$ be the the simplicial map sending all vertices of $K(\G)$, whose colors belong to the same partition class, to one vertex of the standard $2$-simplex; then $H_{b}$ (resp. $H_{r}$) (resp. $H_{g}$) is the preimage by $\mu$ of the star of $v_b$ (resp. $v_{r}$) (resp. $v_{g}$) in $\sigma^\prime.$ 
Therefore it is easy to see that $H_{b}$ is a regular neighbourhood of the (unique) $4$-colored vertex of $K(\G)$ and precisely it is the cone over its disjoint link. 
On the other hand $H_{r}$ (resp. $H_{g}$) is a regular neighbourhood of the $1$-dimensional subcomplex $K_{02}(\G)$ (resp. $K_{13}(\G)$) of $K(\G)$ generated by its 
$i$-colored vertices, with $i\in\{0,2\}$ (resp.  $i\in\{1,3\}$) and it is not difficult to see that $H_{r}$ (resp. $H_{g}$) is a $4$-dimensional handlebody of 
genus $g_{1,3,4} - 1 = m$  (resp. $g_{0,2,4} - 1 = m$), where $m = rk(\pi_1(M^4))$.       

The bijection between $M^4$ and $\widehat{M^4}$ described in Remark  \ref{correspondence-sing-boundary} allows to prove (see \cite{Casali-Cristofori_trisections} for details) that $(H_{b},H_{r},H_{g})$ defines a triple  
$\mathcal T(\G,\varepsilon) = (H_{0},H_{1},H_{2})$ satisfying  Theorem \ref{theorem_B-trisection},  by setting $H_{1}=H_{r}$, $H_2=H_{g}$ and $H_0=H_b$, if $\partial M^4=\emptyset$, or $H_0$ to be a collar of $\partial M^4$ obtained by removing from $H_b$ a suitable neighbourhood of the singular vertex,  if $\partial M^4 \ne \emptyset$.

In fact, $H_{rb}= H_{r} \cap H_{b}$ (resp. $H_{gb}= H_{g} \cap H_{b}$) is the preimage under $\mu$ of the edge of $\sigma^\prime$ depicted in Figure \ref{figsigma} as the ``green'' (resp. ``red'') edge, and turns out to be always an handlebody. 
Moreover, both in the closed and connected boundary case, the central surface $\Sigma= H_{0}\cap H_{1}\cap H_{2} = H_{b}\cap H_{r}\cap H_{g} $ is proved to be a closed connected surface of genus $\rho_{\varepsilon_{\hat 4}}(\Gamma_{\hat{\varepsilon_4}}).$
With regard to $H_{rg}$, this complex is the preimage under $\mu$ of the edge of $\sigma^\prime$ depicted in Figure \ref{figsigma} as the ``blue'' edge.  

By Definition \ref{def_gem-induced trisection}, the triple $\mathcal T(\G,\varepsilon) = (H_{0},H_{1},H_{2})$ of $M^4$ is a gem-induced trisection of $M^4$ if $H_{rg}$  collapses to a graph.

\medskip 

The following proposition states, for weak semi-simple crystallizations, a ``minimality property" regarding the genus of the associated central surface, 
which is actually proved in \cite[Proposition 14]{Casali-Cristofori_trisections} for a larger class of 5-colored graphs. 

\begin{proposition}{\rm \cite{Casali-Cristofori_trisections}} \label{trisection-weak} 
Let $M^4$ be a  compact $4$-manifold with empty or connected boundary and let $\G$ be a weak semi-simple crystallization of $M^4$ with respect to $\varepsilon\in\mathcal P_4.$
Let $\mathcal T (\G,\varepsilon)=(H_{0},H_{1},H_{2})$ be the triple of submanifolds of $M^4$ associated to $\G$ and $\varepsilon$; then  
$$genus(\Sigma(\mathcal T (\G,\varepsilon)))\leq genus(\Sigma(\mathcal T (\bar\G,\bar\varepsilon)))$$ for all $\bar\G\in G_s^{(4)}$ such that $|K(\G)|\cong|K(\bar\G)|$ and for all $\bar\varepsilon\in \mathcal P_4.$
\end{proposition}

It is pointed out in \cite{Chu-Tillman} that, if $(H_0,H_1,H_2)$ is a trisection of a closed $4$-manifold $M^4$ with central surface $\Sigma$, then 
$g(\Sigma)\geq\beta_1(M^4)+\beta_2(M^4).$  

A more general formula (\cite[Proposition 18]{Casali-Cristofori_trisections}) extends the above to a wider class of compact $4$-manifolds with empty or connected boundary. 
By applying it to the case of weak semi-simple crystallizations, and by making use also of Proposition  \ref{weaksemisimple}, we have:

\begin{proposition}{\rm \cite{Casali-Cristofori_trisections}} \label{minimality weak semi-simple}
\ \ 
Let $M^4$ be a compact $4$-manifold with empty or connected boundary which admits a weak semi-simple crystallization $\G$ giving rise to a gem-induced trisection $\mathcal T(\G,\varepsilon)$.
Then, 
$$g_{GT}(M^4)\ =\   \frac 1 2 ( \rho_{\varepsilon}(\G) + m) \ =\  \beta_2(M^4)+\beta_1(M^4)+ 2\left(m-\beta_1(M^4)\right),$$ 
\noindent with $m = rk(\pi_1(M^4)).$

\noindent In particular, $$g_{GT}(M^4) \ = \ \frac 12\rho_\varepsilon(\G) \ = \ \beta_2(M^4) $$ for each compact (simply-connected) $4$-manifold $M^4$, with empty or connected boundary, which admits a weak simple crystallization $\G$ giving rise to a gem-induced trisection $\mathcal T(\G,\varepsilon)$.

In this case, if $M^4$ is closed, its G-trisection genus (equal to the second Betti number $\beta_2(M^4)$) coincides with its trisection genus.  
\end{proposition}

Finally, given a $5$-colored graph $\G\in G_s^{(4)}$ representing a compact $4$-manifold $M^4$ with empty or connected boundary, a sufficient condition  is known for $\mathcal T (\G,\varepsilon)$ to be a B-trisection of $M^4$ for each cyclic permutation $\varepsilon \in \mathcal P_4$. It makes use of the existence of a presentation of $\pi_1 (\widehat M^4)$  with generator set in bijection with 4-colored edges of $\Gamma$  and relator set in bijection with bicolored cycles of $\G$ involving color $4$: see \cite{Casali-Cristofori_generalized} for details.

\begin{proposition}{\rm \cite{Casali-Cristofori_trisections}}  \label{CS gem-induced trisections}
Let $< X,R>$ be the presentation of $\pi_1 (\widehat M^4)$  with $X=\{x_1, \dots, x_p\}$ in bijection with 4-colored edges of $\Gamma$  and $R=\{r_1, \dots, r_q\}$ in bijection with $\{4,i\}$-cycles of $\G$, for each $i \in \Delta_3$. 

If the presentation  $< X,R>$ collapses to the trivial one through a finite sequence of moves of the following type: 
\begin{equation*} \label{CNS}
 \begin{aligned} \text{if} \ \  r_j=x_s \ (j \in\mathbb N_q, \  s \in\mathbb N_p), \ \  & \text{then delete} \ x_s \ \text{from the generator set,} \\
\ & \text{and from each relation containing it, too, } 
\end{aligned} 
\end{equation*}
\noindent 
then $M^4$ admits a gem-induced trisection $\mathcal T(\Gamma, \varepsilon)$, for each $\varepsilon \in \mathcal P_4$.  
\end{proposition}
\bigskip

All simple, semi-simple or weak semi-simple crystallizations depicted in Figures \ref{fig.S4}--\ref{fig.S2xS2} and 
\ref{fig.S2xD2}--\ref{fig.csi_c} turn out to satisfy the sufficient condition of Proposition \ref{CS gem-induced trisections}; 
the same happens with the simple crystallization of the surface K3 presented in \cite{Basak-Spreer}. 
Therefore, all of them give rise to gem-induced trisections, that, by Proposition \ref{minimality weak semi-simple}, realize the G-trisection genus of the represented manifolds (see 
Tables 2, where only manifolds $M^4$ such that $\pi_1(\widehat{M^4})=0$ are taken into account). 

Proposition \ref{minimality weak semi-simple} also ensures that, for the closed simply-connected $4$-manifolds of Tables 1 and 2, 
the described crystallizations turn out to realize also the trisection genus.

\begin{table}[!h]

\begin{center}
\begin{tabular}{|c||c|} 
\hline \ & \ \cr
$M^4$ & $g_{GT}(M^4)$ \cr 
\hline \hline  
$\mathbb S^4$ & 0 \cr  \hline
$\mathbb {CP}^2$ & 1 \cr  \hline
$\mathbb{S}^{2} \times \mathbb{S}^{2}$ & 2 \cr  \hline
$K3$ & 22 \cr  \hline
$\mathbb{S}^{2} \times \mathbb{D}^{2}$ & 1 \cr  \hline
$\mathbb Y^4_h$ & $h$ \cr  \hline
$\xi_c$ ($c \ge 2$)   
& 1 \cr  \hline
\end{tabular}
\caption{computation of the G-trisection genus.} 
\end{center}
\label{table_2}
\end{table}

\smallskip

The following proposition shows that the G-trisection genus has the same behaviour as the regular genus  with respect to connected sums (compare with Proposition \ref{prop:connected_sum}):

\begin{proposition} {\em \cite{Casali-Cristofori_trisections}} \label{g_GT properties} \ 
Let $M^4_1, M^4_2$ be two compact $4$-manifolds with empty or connected boundary admitting gem-induced trisections. Then $M^4_1 \sharp M^4_2$ admits gem-induced trisections, too, and 
$$\ \ g_{GT}(M^4_1 \sharp M^4_2) \le g_{GT}(M^4_1) + g_{GT}(M^4_2).$$ 

\noindent Furthermore, equality holds if  $M^4_1$ and $M^4_2$ admit B-trisections induced by weak semi-simple crystallizations.
\end{proposition} 

\begin{remark} {\em As a consequence of the previous proposition, for any compact $4$-manifold $M^4$, with empty or connected boundary, that is a connected sum of the manifolds in Table 2,  
the equality $$g_{GT}(M^4)=\beta_2(M^4)+\beta_1(M^4)$$ holds. 
\noindent In particular, for all closed simply-connected ``standard'' $4$-manifolds, Propositions \ref{minimality weak semi-simple} and \ref{g_GT properties} ensure that the trisection genus equals the second Betti number, as proved by Spreer and Tillmann in \cite{Spreer-Tillmann}
by using simple crystallizations.}
\end{remark}

\bigskip

\nocite{*}
\bibliographystyle{plain}
\bibliography{weak-semisimple_arxiv}

\begin{thebibliography}{10}

\bibitem{Bandieri-Cristofori-Gagliardi}
P.~Bandieri, P.~Cristofori, and C.~Gagliardi.
\newblock Nonorientable 3-manifolds admitting coloured triangulations with at
  most 30 tetrahedra.
\newblock {\em Journal of Knot Theory and its Ramifications}, 18:381--395,
  2009.

\bibitem{Basak}
B.~Basak.
\newblock Genus-minimal crystallizations of {PL} 4-manifolds.
\newblock {\em Beitr. Algebra Geom.}, 59(1):101--111, 2018.

\bibitem{Basak-Casali}
B.~Basak and M.R. Casali.
\newblock Lower bounds for regular genus and gem-complexity of {PL}
  4-manifolds.
\newblock {\em Forum Math.}, 29(4):761--773, 2017.

\bibitem{Basak-Spreer}
B.~Basak and J.~Spreer.
\newblock Simple crystallizations of 4-manifolds.
\newblock {\em Adv. Geom.}, 16(1):111--130, 2016.

\bibitem{Bell-et-al}
M.~Bell, J.~Hass, J.H. Rubinstein, and S.~Tillmann.
\newblock Computing trisections of 4-manifolds.
\newblock {\em Proc. Nat. Acad. Sci. USA}, 115(43):10901--10907, 2018.

\bibitem{Casali_Forum2003}
M.R. Casali.
\newblock On the regular genus of 5-manifolds with free fundamental group.
\newblock {\em Forum math.}, 15:465--475, 2003.

\bibitem{Casali-Cristofori_JKTR2008}
M.R. Casali and P.~Cristofori.
\newblock A catalogue of orientable 3-manifolds triangulated by $30$ coloured
  tetrahedra.
\newblock {\em J. Knot Theory Ramifications}, 17:1--23, 2008.

\bibitem{Casali-Cristofori_CATALOGUE}
M.R. Casali and P.~Cristofori.
\newblock Cataloguing {PL} 4-manifolds by gem-complexity.
\newblock {\em Electron. J. Combin.}, 22(4):\#P4.25, 1--25, 2015.

\bibitem{Casali-Cristofori_trisections}
M.R. Casali and P.~Cristofori.
\newblock {\it Gem-induced trisections of compact {PL} 4-manifolds}, preprint
  2019 (arXiv:1910.08777).

\bibitem{Casali-Cristofori_generalized}
M.R. Casali and P.~Cristofori.
\newblock {\it Classifying compact 4-manifolds via generalized regular genus
  and {G}-degree}, preprint 2019 (arXiv:1912.01302).

\bibitem{CC_handledec-weaksimple}
M.R. Casali and P.~Cristofori.
\newblock {\it Compact 4-manifolds admitting special handle-decompositions},
  preprint 2020.

\bibitem{Casali-Cristofori-Dartois-Grasselli}
M.R. Casali, P.~Cristofori, S.~Dartois, and L.~Grasselli.
\newblock Topology in colored tensor models via crystallization theory.
\newblock {\em J. Geom. Phys.}, 129:142--167, 2018.

\bibitem{Casali-Cristofori-Gagliardi_MONTESINOS}
M.R. Casali, P.~Cristofori, and C.~Gagliardi.
\newblock {\em Classifying {PL} 4-manifolds via crystallizations: results and
  open problems}.
\newblock in: ``Mathematical Tribute to Professor Jos\'e Maria Montesinos
  Amilibia''. Universidad Complutense Madrid, [{ISBN}: 978-84-608-1684-3]
  edition, 2016.

\bibitem{Casali-Cristofori-GagliardiJKTR2015}
M.R. Casali, P.~Cristofori, and C.~Gagliardi.
\newblock {PL} 4-manifolds admitting simple crystallizations: framed links and
  regular genus.
\newblock {\em Journal of Knot Theory and its Ramifications}, 25(1):165005 (14
  pages), 2016.

\bibitem{Casali-Cristofori-Grasselli}
M.R. Casali, P.~Cristofori, and L.~Grasselli.
\newblock G-degree for singular manifolds.
\newblock {\em RACSAM}, 112(3):693--704, 2018.

\bibitem{Casali-Gagliardi_ProcAMS}
M.R. Casali and C.~Gagliardi.
\newblock Classifying {PL} 5-manifolds up to regular genus seven.
\newblock {\em Proc. Amer. Math. Soc.}, 120:275--283, 1994.

\bibitem{Casali-Grasselli2019}
M.R. Casali and L.~Grasselli.
\newblock Combinatorial properties of the {G}-degree.
\newblock {\em Rev. Mat. Complut.}, 32(1):239--254, 2019.

\bibitem{Chu-Tillman}
M.~Chu and S.~Tillmann.
\newblock Reflections on trisection genus.
\newblock {\em Rev. Roumaine Math. Pures Appl. (Proceedings of JARCS 2017)},
  64(4), 2019.

\bibitem{Cristofori-Fomynikh-Mulazzani-Tarkaev}
P.~Cristofori, E.~Fomynikh, M.~Mulazzani, and V.~Tarkaev.
\newblock 4-colored graphs and knot/link complements.
\newblock {\em Results Math.}, 72(1-2):471--490, 2017.

\bibitem{Cristofori-Mulazzani}
P.~Cristofori and M.~Mulazzani.
\newblock Compact 3-manifolds via 4-colored graphs.
\newblock {\em RACSAM}, 110(2):395--416, 2015.

\bibitem{Ferri-Gagliardi}
M.~Ferri and C.~Gagliardi.
\newblock The only genus zero $n$-manifold is ${S}^n$.
\newblock {\em Proc. Amer. Math. Soc.}, 85:638--642, 1982.

\bibitem{FGG}
M.~Ferri, C.~Gagliardi, and L.~Grasselli.
\newblock A graph-theoretical representation of {PL}-manifolds. a survey on
  crystallizations.
\newblock {\em Aequationes Math.}, 31:121--141, 1986.

\bibitem{Gagliardi_1981}
C.~Gagliardi.
\newblock Extending the concept of genus to dimension $n$.
\newblock {\em Proc. Amer. Math. Soc.}, 81:473--481, 1981.

\bibitem{Gay-Kirby}
D.~Gay and R.~Kirby.
\newblock Trisecting 4-manifolds.
\newblock {\em Geom. Topol.}, 20:3097--3132, 2016.

\bibitem{Gurau-book}
R.~Gurau.
\newblock {\em Random {T}ensors}.
\newblock Oxford University Press, 2016.

\bibitem{Lins}
S.~Lins.
\newblock {\em Gems, computers and attractors for 3-manifolds.}
\newblock Number~5 in Knots and Everything. World Scientific, 1995.

\bibitem{Mandelbaum}
R.~Mandelbaum.
\newblock Four-dimensional topology: an introduction.
\newblock {\em Bull. Amer. Math. Soc.}, 2:1--159, 1980.

\bibitem{Pezzana_1}
M.~Pezzana.
\newblock Sulla struttura topologica delle variet\'a compatte.
\newblock {\em Atti Semin. Mat. Fis. Univ. Modena}, 23:269--277, 1974.

\bibitem{Pezzana_2}
M.~Pezzana.
\newblock Diagrammi di {H}eegaard e triangolazione contratta.
\newblock {\em Boll. Un. Mat. Ital.}, 12:98--105, 1975.

\bibitem{Rubinstein-Tillmann}
H.~Rubinstein and S.~Tillmann.
\newblock {\it Multisections of piecewise linear manifolds}, to appear in
  Indiana Univ. Math. J.

\bibitem{Spreer-Tillmann}
J.~Spreer and S.~Tillmann.
\newblock The trisection genus of standard simply connected {PL} 4-manifolds.
\newblock In in: Leibniz International Proceedings~in Informatics~(LIPICS),
  editor, {\em 34nd International Symposium on Computational Geometry (SoCG
  2018)}, volume~99, pages 1--71, 2018.

\end{thebibliography}
\end{document}